\documentclass[15pt, twoside]{article}
\usepackage{amssymb}
\usepackage[all]{xy}
\usepackage{amsthm,amsmath,color}
\usepackage{pgfmath}
\usepackage{tikz-cd}
\usetikzlibrary{arrows, matrix}
\usepackage{zref-abspage}
\usepackage{perpage}
\usepackage[top=2.5cm,right=3.5cm,bottom=2.5cm,left=2.5cm]{geometry}
\usepackage{graphicx}
\usepackage{fancyhdr}
\usepackage{makeidx}
\usepackage{ifpdf}
\usepackage[colorlinks=true, linkcolor=blue, linktoc=all, citecolor=magenta, filecolor=cyan, urlcolor=cyan, linkbordercolor={1 0 0}, citebordercolor={0 1 0}, urlbordercolor={0 1 1}]{hyperref}
\newtheorem{theorem}{Theorem}[section]
\newtheorem{definition}{Definition}[section]
\newtheorem{proposition}{Proposition}[section]
\newtheorem{lemma}{Lemma}[section]
\newtheorem{corollary}{Corollary}[section]
\theoremstyle{definition}
\newtheorem{example}{Example}[section]
\begin{document}
\renewcommand{\headrulewidth}{1pt}
\fancyhead{}
\fancyfoot{}
\chead[M. Roshandelbana, S. Varsaie]{Analytic Approach To $\nu$-classes}
\lhead[\thepage]{}
\rhead[]{\thepage}
\pagestyle{fancy}
\author{M. Roshandelbana, S. Varsaie
\\
\small{
\textit{Department of Mathematics,}}
\\
\small{\textit{Institute for Advanced Studies in Basic Sciences (IASBS),}}
\\
\small{\textit{Zanjan 45137-66731, Iran.}}}
\title{Analytic Approach To a Generalization of Chern Classes in
Supergeometry}
\renewcommand{\baselinestretch}{1.5}
\setlength{\baselineskip}{1.5\baselineskip}
\setcounter{page}{1}
\maketitle
\tableofcontents
\begin{abstract}
Some cohomology elements, called $\nu$-classes, as a supergeneralization of universal Chern classes, are introduced for canonical super line bundles over $\nu$-projective spaces, a novel supergeometric generalization of projective spaces. It is shown that these classes may be described by analytic representatives of elements of generalized de Rham cohomology.
\\
\textbf{Keywords:}
Chern classes; Super projective spaces; $\nu$-projective spaces; Canonical super line bundles.
\end{abstract}
\section{Introduction}
In this paper, we introduce an analytic approach to study the first $\nu$-classes. These classes are a novel generalization of the first Chern classes in supergeometry introduced in \cite{Dr2}. In addition, at the end, we introduce some cocycles in a de Rham complex which may be considered as representatives of $\nu$-classes of rank higher than 1 for universal super vector bundles.
\\
In fact, among geometric concepts, Chern classes are important topological invariants associated to complex vector bundles. So it is important to obtain a proper generalization of Chern classes in supergeometry. Since 1980, many efforts have been made to generalize this concept which we now mention some of them.
\\
In \cite{Bartocci}, \cite{Bartocci book}, \cite{Bruzzo} and \cite{Landi}, supermanifolds are defined in the sense of De Witt \cite{DeWitt} and Rogers \cite{Rogers}. Note that in this approach, a super vector bundle over a De Witt supermanifold does always carry a connection. In all of these references, for a given complex super vector bundle of rank $(r, s)$ over a De Witt supermanifold, $r$ even and $s$ odd Chern classes are defined and consequently, one may have even and odd total Chern classes. In addition, a representative of Chern classes in terms of the curvature form of a connection is given. Note that the second author has shown, in \cite{Dr5}, that the even and odd Chern classes of a super vector bundle are the same as Chern classes of two canonically associated vector bundles. Since a category with pairs of vector bundles as its objects is not equivalent to the category of super vector bundles, this theory of Chern classes for super vector bundles is exactly the common one for the former category. For more details, see \cite{Dr5}.
\\
In the paper \cite{Manin and Penkov}, based on the Leites \cite{Leites} and Kostant \cite{Kostant} approach to supermanifolds, a left(right) connection is introduced and its associated curvature form is computed. Then, considering a locally free sheaf with a smooth left(right) connection and applying invariant polynomials to its curvature form, a globally defined closed differential form is obtained. It is shown that there is a representative of the Chern classes in terms of the introduced differential forms. But as it is mentioned in \cite{Manin and Penkov}, the classes obtained in this way are nothing but the Chern classes of the reduced vector bundle. Indeed, in the smooth construction, the category of locally free sheaves on a supermanifold is equivalent to the category of $\mathbb{Z}_{2}$-graded locally free sheaves on the reduced manifold. So the introduced Chern classes are just those of common geometry.
\\
In \cite{Voronov-Manin-Penkov}, associated to a super vector bundle $\mathcal{E}$, it is considered a decomposition as $\mathcal{E}_{red} = \mathcal{E}_{0} \oplus \mathcal{E}_{1}$ for the reduced vector bundle. Then, by an approach which coincides to that of Quillen(\cite{Quillen}), a representative of Chern classes of $\mathcal{E}$ is given in terms of the Chern classes of $\mathcal{E}_{0}$ and $\mathcal{E}_{1}$. Since $\mathcal{E}_{0}$ and $\mathcal{E}_{1}$ are components of the reduced vector bundle, this representative do not have any information about the superstructure, i.e., this theory does not indicate the extent to which even and odd elements of a super vector bundle are associated to each other. Solving this problem is a motivation for defining $\nu$-classes.
\\
This paper aims at describing $\nu$-classes by analytic representatives of elements of generalized de Rham cohomology. The idea of analytic approach to describe $\nu$-classes derives from the idea of analytic description of Chern classes in common geometry where there are different ways to approach Chern classes, such as homotopy or analytic approach. In homotopy approach, one can consider any vector bundle $\xi$ on a manifold $M$ as a pullback of the universal bundle along a map, say $f$, from $M$ to infinite Grassmannian, unique up to homotopy. Therefore, one can describe Chern classes of $\xi$ as the pullback of Chern classes of the universal bundle. In \cite{Afshari}, this theory is generalized in supergeometry. Indeed, a canonical super vector bundle over $\nu$-grassmannian is constructed and it is shown that it has the property of the universal super vector bundle. So by having universal $\nu$-classes, one can obtain $\nu$-classes of rank higher than $1$ for any super vector bundle. 
\section{Preliminaries}
In this section, we deal with some basic definitions of supergeometry such as supermanifolds, $\nu$-domains, derivations and etc.
\subsection{Supermanifolds}
A \textit{super ringed space} is a pair $(M, \mathcal{O})$ where $M$ is a topological space and $\mathcal{O}$ is a sheaf of supercommutative $\mathbb{Z}_{2}$-graded rings. So for any open subset $U \subset M$, the section $\mathcal{O}(U)$ is a supercommutative $\mathbb{Z}_{2}$-graded ring. We denote its even and odd parts by $\mathcal{O}^{0}(U)$ and $\mathcal{O}^{1}(U)$ respectively. Thus, one has
\begin{equation*}
\mathcal{O}(U) = \mathcal{O}^{0}(U) \oplus \mathcal{O}^{1}(U) 
\end{equation*}
with the property that $\mathcal{O}^{i}.\mathcal{O}^{j} \subset \mathcal{O}^{(i+j) \bmod 2}$.
\\
We call each element of  $\mathcal{O}^{0}$ and $\mathcal{O}^{1}$ a homogeneous element of even and odd parity respectively.
\\
We denote the parity of a homogeneous element $a$ by $p(a)$ which equals $0$ if it is even and equals $1$ if it is odd.
\\
Now, let $(M, \mathcal{O}_{M})$ and $(N, \mathcal{O}_{N})$ be two super ringed spaces. By a morphism of super ringed spaces from $(M, \mathcal{O}_{M})$ to $(N, \mathcal{O}_{N})$, we mean a pair $\psi = (\widetilde{\psi}, \psi^{*})$ where $\widetilde{\psi}: M \longrightarrow N$ is a continuous map and $\psi^{*}: \mathcal{O}_{N} \longrightarrow \widetilde{\psi}_{*}(\mathcal{O}_{M})$ is a homomorphism between sheaves of $\mathbb{Z}_{2}$-graded rings. Thus, for any open subset $V \subset N$, we have
\begin{equation*}
\psi^{*}_{V}: \mathcal{O}_{N}(V) \longrightarrow \mathcal{O}_{M}\big({\widetilde{\psi}}^{-1}(V)\big)
\end{equation*}
which commutes with the corresponding restriction morphisms and in addition satisfies the following:
\begin{equation*}
\psi^{*}_{V}(\mathcal{O}^{0}_{N}(V)) \subset \mathcal{O}^{0}_{M}(\widetilde{\psi}^{-1}(V)), \qquad \psi^{*}_{V}(\mathcal{O}^{1}_{N}(V)) \subset \mathcal{O}^{1}_{M}(\widetilde{\psi}^{-1}(V)).
\end{equation*}
Let $U$ be an open subset of $\mathbb{R}^{m}$. A \textit{superdomain} $U^{m|n}$ is a super ringed space $(U, \mathcal{C}^{m|n}_{U})$ where $\mathcal{C}^{m|n}_{U}$ is the sheaf of supercommutative rings such that for any open subset $V \subset U$ we have
\begin{equation}\label{2}
\mathcal{C}^{m|n}_{U}(V) = \mathbf{C}^{\infty}(V)[e_{1}, \cdots, e_{n}]
\end{equation}
where $\{e_{1}, \cdots, e_{n}\}$ are anticommuting variables satisfying $e_{i}e_{j} = -e_{j}e_{i}$ and $\mathbf{C}^{\infty}(V)$ is the algebra of smooth functions on $V$.
\\
If there is no ambiguity, we write $\mathcal{C}_{U}$ for $\mathcal{C}_{U}^{m|n}$.
\\
From now on, let $F= \mathbb{R}$ or $\mathbb{C}$. Suppose $\{z_{k}\}_{1 \leq k \leq m}$ be a global coordinate system on $F^{m}$. Then, it is called even coordinates on $U^{m|n}$. In addition, $\{e_{l}\}$ introduced by \eqref{2} is called odd coordinates on $U^{m|n}$.
\\
A super ringed space $(M, \mathcal{O})$ is called a \textit{supermanifold} of dimension $m|n$ if it is locally isomorphic to $F^{m|n}=(F^{m}, \mathcal{C}_{F^{m}})$, namely, there exists an open cover $\{V_{\alpha}\}$ for $M$ such that for any $\alpha$, there exists an isomorphism $\psi_{\alpha} =(\widetilde{\psi}_{\alpha}, \psi^{*}_{\alpha})$ from $(V_{\alpha}, \mathcal{O}_{|_{V_{\alpha}}})$ to $F^{m|n}$. Thus for an open subset $V \subset F^{m}$, $\widetilde{\psi}_{\alpha}: V_{\alpha} \longrightarrow F^{m}$ is a homeomorphism and $\psi^{*}_{\alpha}:\mathcal{C}_{F^{m}}(V) \longrightarrow \mathcal{O}\big(\widetilde{\psi}^{-1}_{\alpha}(V)\big)$ is an isomorphism of $\mathbb{Z}_{2}$-graded rings.\\
Clearly, the morphisms between supermanifolds are just morphisms between the corresponding super ringed spaces such that for any $x \in V_{\alpha}$, $\psi^{*}: \mathcal{C}_{F^{m}, \widetilde{\psi}(x)} \to \mathcal{O}_{x}$ is local, i.e., $\psi^{*}(m_{\widetilde{\psi}(x)}) \subset m_{x}$, where $m_{x}$ is the unique maximal ideal in $\mathcal{O}_{x}$.
\\
A \textit{super vector bundle} of rank $r|s$ over a supermanifold $(M, \mathcal{O})$ is a locally free sheaf of $\mathcal{O}$-modules of rank $r|s$ over $M$.
\subsection{$\nu$-domains}
By a $\nu$\textit{-domain} $(U, \mathcal{C}_{U}, \nu)$, we mean a superdomain $(U, \mathcal{C}_{U})$ with an odd involution $\nu: \mathcal{C}_{U} \longrightarrow \mathcal{C}_{U}$ when it is considered as a morphism between sheaves of $\mathbf{C}^{\infty}_{U}$-modules. In other words, one has
\begin{equation*}
\nu^{2} = 1, \qquad \nu (\mathcal{C}_{U})^{0} \subset \big(\mathcal{C}_{U}\big)^{1}, \qquad \nu \big(\mathcal{C}_{U}\big)^{1} \subset \big(\mathcal{C}_{U}\big)^{0}.
\end{equation*}
Hence,  $\mathcal{C}_{U}$ may be considered as a sheaf of $\mathbb{C}_{\nu_{0}}$-modules where $\mathbb{C}_{\nu_{0}} = \mathbb{C}[\nu_{0}]$ is a ring generated by indeterminate $\nu_{0}$ with $\nu_{0}^{2}=1$. Indeed, one may define $\nu_{0} a:= \nu(a)$, for each $a \in \mathcal{C}_U$.
\\
A morphism of $\nu$-domains is a morphism, say $\psi=(\widetilde{\psi}, \psi^{*})$, from $(F^{m}, \mathcal{C}_{F^{m}}^{m|n}, \nu)$ to $(F^{k}, \mathcal{C}_{F^{k}}^{k|l}, \nu^{\prime})$ where $\psi^{*}$ preserves the $\mathbb{C}_{\nu_{0}}$-module structure of corresponding sheaves, i.e., for any open subset $V \subset F^{k}$, the morphism $\psi^{*}_{V}: \mathcal{C}_{F^{k}}^{k|l}(V) \longrightarrow  \mathcal{C}_{F^{m}}^{m|n}\big(\widetilde{\psi}^{-1}(V)\big) $ is a $\mathbb{C}_{\nu_{0}}$-module homomorphism.
\subsection{Differential forms}
Now, we introduce the notion of derivations and differential forms in supergeometry.
\\
A \textit{homogeneous derivation} of a superalgebra $\mathcal{A}$ on a field $k$ is a $k$-linear map $D: \mathcal{A} \longrightarrow \mathcal{A}$ such that 
\begin{equation*}
D(ab) = (Da)b + (-1)^{p(D)p(a)}a(Db) \qquad a, b \in \mathcal{A},
\end{equation*}
where by $p(D)$ we mean the parity of the derivation.
\\
A map $D: \mathcal{A} \longrightarrow \mathcal{A}$ is a \textit{derivation} of $\mathcal{A}$ if any of its homogeneous components is a homogeneous derivation. The space of all derivations on $\mathcal{A}$ is denoted by $Der(\mathcal{A})$.
\\
Let $(M, \mathcal{O})$ be a supermanifold. In addition, let $Der(\mathcal{O}(U))$ be the space of all derivations of $\mathcal{O}(U)$ for an open subset $U \subset M$. Considering $T(U)$ as the tensor algebra of $Der(\mathcal{O}(U))$ over $\mathcal{O}(U)$ and $T^{r}(U)$ as the homogeneous subspace corresponding to $r \in \mathbb{N}$, the set $Hom _{\mathcal{O}(U)}(T^{r}(U), \mathcal{O}(U))$ contains of all $r$-linear maps on $Der(\mathcal{O}(U))$ such that for any $X_{1}, \cdots, X_{r} \in Der(\mathcal{O}(U))$, $f \in \mathcal{O}(U)$ and $\omega \in Hom _{\mathcal{O}(U)}(T^{r}(U), \mathcal{O}(U))$, we have
\begin{equation*}
\big<X_{_{1}}, \cdots, fX_{_{l}}, \cdots, X_{_{r}} | \omega\big> = (-1)^{p(f) \sum_{i=1}^{l-1} p(X_{i})} f \big<X_{_{1}}, \cdots, X_{_{l}}, \cdots, X_{_{r}}| \omega\big>.
\end{equation*}
\begin{definition}
The set of $r$\textit{-differential forms} on $U$ denoted by $\Omega^{r}(U, \mathcal{O})$ is the set of all maps $\omega \in Hom _{\mathcal{O}(U)}(T^{r}(U), \mathcal{O}(U))$ such that it satisfies the following condition:
\begin{equation*}
\big< X_{_{1}}, \cdots, X_{_{l}}, X_{_{l+1}}, \cdots, X_{_{r}}| \omega\big> = (-1)^{1+p(X_{l})p(X_{l+1})} \big< X_{_{1}}, \cdots, X_{_{l+1}}, X_{_{l}}, \cdots, X_{_{r}} | \omega\big>.
\end{equation*}
\end{definition}
$\Omega^{r}(U, \mathcal{O})$ is an $\mathcal{O}(U)$-module with the property below
\begin{equation*}
\big< X_{_{1}}, \cdots, X_{_{r}} | \omega f\big>=\big< X_{_{1}}, \cdots, X_{_{r}}| \omega \big> f,
\end{equation*}
In addition, $\Omega^{r}(U, \mathcal{O})$ is $\mathbb{Z}_{2}$-graded. Thus, we have 
\begin{equation*}
\big< X_{_{1}}, \cdots, X_{_{r}} | \omega \big> \in \mathcal{O}^{k}(U),
\end{equation*}
where $k:= \big(p(\omega) + \sum_{i=1}^{r} p(X_{i})\big) \bmod 2$.
\\
Set
\begin{equation*}
\Omega(U, \mathcal{O})= \oplus_{r=0}^{\infty}\Omega^{r}(U, \mathcal{O}), \qquad \Omega^{0}(U, \mathcal{O}) = \mathcal{O}(U)
\end{equation*}
The correspondence $U \longrightarrow \Omega(U, \mathcal{O})$ defines a sheaf of supercommutative superalgebras on $M$ called the sheaf of differential forms.
\begin{lemma}\label{derivation}
Let $U \subset M$ be an arbitrary open set. There exists a unique even derivation $d: \Omega(U, \mathcal{O}) \longrightarrow \Omega(U, \mathcal{O})$ of degree $1$, satisfying the following properties:
\begin{itemize}
\item[1)]
On $\Omega^{0}(U, \mathcal{O})$, we have
\begin{align*}
&d: \Omega^{0}(U, \mathcal{O}) \longrightarrow \Omega^{1}(U, \mathcal{O})
\\
&\big< X | dg\big> = X g, \qquad \forall g \in \mathcal{O}(U)=\Omega^{0}(U, \mathcal{O}), \quad X \in Der(\mathcal{O}(U)).
\end{align*}
\item[2)]
$d^{2}= 0$.
\end{itemize}
\end{lemma}
\section{$\nu$-projective spaces}\label{section}
In this section, first, we introduce a new generalization of projective spaces in supergeometry, different from superprojective spaces, called $\nu$-projective spaces which are constructed by gluing $\nu$-domains. At the end, we introduce a $1|0$-super line bundle over this space.
\subsection{Construction of \texorpdfstring{$\nu$}{Lg}- projective spaces}
In this section, by introducing gluing morphisms, it is shown the way a $\nu$-projective space is constructed by gluing a number of copies of $\nu$-domains $(\mathbb{C}^{m}, \mathcal{C}_{\mathbb{C}^{m}}, \nu)$, where $\mathcal{C}_{\mathbb{C}^m}$ is the sheaf of complex valued smooth functions.
\\
For each $i \in \{1, \cdots, m+n+1\}$, let $(U_{i}, \mathcal{O}_{i}) = (\mathbb{C}^{m}, \mathcal{C}_{\mathbb{C}^{m}})$ be a $\nu$-domain. We write $U_{i}$ instead of $(U_{i}, \mathcal{O}_{i})$ for brevity.
\\
The $\nu$-domain $U_{i}$ is called \textit{standard} if $1 \leq i \leq m +1$, otherwise it is called \textit{nonstandard}. The reason for the adoption of this convention will be further clarified.
\\
Corresponding to the index $i$, we label the $\nu$-domain $U_{i}$ by an even supermatrix, say $A_{i}$, of dimension $1|0 \times (m+1)|n$ with the property that the $i$-th entry is equal to $1$ if $1 \leq i \leq m+1$, otherwise it is equal to $1\nu$ which is a formal symbol. Obviously, $A_{i}$ may be decomposed into two blocks $B_{1}$ and $B_{2}$ which are respectively called even and odd part of $A_{i}$ because of the parity of their elements. In addition, except the $i$-th entry the blocks $B_j$, ($j=1,2$), are filled by even  and odd coordinates from left to right respectively according to the order induced by their indices. In this process, a coordinate, say $w$, is replaced by $\nu(w)$ if it places in a part with opposite parity. We separate even and odd parts of $A_{i}$ by a vertical line called \textit{divider line}.
\\
Denote the $j$-th entry of $A_{i}$ by $M_{j}(A_{i})$. It is clear that $M_{j}(A_{i})$ is even if $1 \leq j \leq m+1$, otherwise it is an odd entry of $A_{i}$. Obviously, $M_{i}(A_{i})=1$ for $i \leq m+1$ and $M_{i}(A_{i})=1\nu$ for $i > m+1$.
\\
Now, we introduce $M^{\prime}_{j}(A_{i})$ which will be used in the subsequent sections. 
\\
For $j \leq m+1$, set $M^{\prime}_{j}(A_{i}) = M_{j}(A_{i})$ and otherwise define $M^{\prime}_{j}(A_{i}) = \nu\big(M_{j}(A_{i})\big)‌$. Clearly, $M^{\prime}_{j}(A_{i})$ always takes even values and $M^{\prime}_{i}(A_{i})$ is always equal to the identity. By convention $\nu^{0}=id$, we may set $M^{\prime}_{j}(A_{i}) = \nu^{p(j)} \big(M_{j}(A_{i})\big)$ for each $j$ where $p(j)=0$ if $j \leq m+1$, otherwise $p(j)=1$.
\subsubsection{Gluing morphisms of $\nu$-domains}
Now, we are going to introduce gluing morphisms of $\nu$-domains.
\\
Consider $\nu$-domains $(U_{i}, \mathcal{O}_{i})$ and $(U_{j}, \mathcal{O}_{j})$ labeled by $A_{i}$ and $A_{j}$ respectively. Now, let $U_{ji}$ be a set consisting of those points of the $\nu$-domain $U_{j}$ for which $M^{\prime}_{i}(A_{j})$ is invertible. The gluing morphism $g_{ij}$ on $U_{ji}$ is as follows:
\begin{equation} \label{tilde}
g_{ij} = (\widetilde{g}_{ij}, g^{*}_{ij}): (U_{ji}, \mathcal{O}_{j}|_{U_{ji}}) \longrightarrow  (U_{ij}, \mathcal{O}_{i}|_{U_{ij}})
\end{equation}
where $\widetilde{g}_{ij}: U_{ji} \longrightarrow U_{ij}$ is a homeomorphism to be defined by lemma \ref{gtilda} and $g^{*}_{ij}: \mathcal{O}_{i}|_{U_{ij}} \longrightarrow \mathcal{O}_{j}|_{U_{ji}}$ is an isomorphism between sheaves determined by defining on each entry of $D_{i}A_{i}$ as a rational expression which appears as the corresponding entry provided by the pasting equation(see \cite{Varadarajan}) :
\begin{equation} \label{DJ}
D_{i}\Big( \big(M^{\prime}_{i}(A_{j})\big)^{-1}A_{j} \Big) = D_{i} A_{i}
\end{equation}
where by $D_{i}(A_{i})$ we mean the supermatrix $A_{i}$ in which the $i$-th entry, $M_{i}(A_{i})$, is omitted.
\\
By \cite{Varadarajan}, the uniqueness of the morphism $g_{ij}$ is concluded.

\begin{proposition} \label{gluing condition}
The morphisms $\{g_{ij}^{*}\}$ defined as above, are gluing morphisms. In other words, the corresponding sheaves of $\nu$-domains $U_{i}$ and $U_{j}$ can be glued together through the morphism $g_{ij}^{*}$.
\end{proposition}
\textit{Proof.}
By (\cite{Varadarajan}, page135), $\{g_{ij}^{*}\}$ are gluing morphisms if and only if the following conditions hold:
\begin{itemize}
\item[1)]
$g^{*}_{ii} = id$
\item[2)]
$g_{ji}^{*} \circ  g_{ij}^{*} = id$
\item[3)]
$g_{ki}^{*} \circ g_{jk}^{*} \circ g_{ij}^{*}= id$
\end{itemize}
\begin{itemize}
\item[1)]
By definition, we know that the morphism $g_{ii}^{*}$ is defined by the equation $D_{i}\Big( \big(M^{\prime}_{i}(A_{i})\big)^{-1}A_{i} \Big) = D_{i} A_{i}$. Since $M^{\prime}_{i}(A_{i}) = 1$, the equality holds.
\item[2)]
Let $i \leq m+1$ and $j > m+1$. The morphism $g^{*}_{ji}$ is obtained from the following:
\begin{equation}\label{equation of Aj}
D_{j}\Big(\big(M^{\prime}_{j}(A_{i})\big)^{-1}A_{i} \Big) = D_{j} A_{j}.
\end{equation}
To compute $g_{ji}^{*} \circ g_{ij}^{*}$, it is required to replace the supermatrix $A_{j}$ in \eqref{DJ} by the left hand side of \eqref{equation of Aj}. So one has
\begin{equation}\label{g}
D_{i}\Big( \big(M^{\prime}_{i} \big[Z^{-1}A_{i}\big]_{j}\big)^{-1}\big[Z^{-1}A_{i}\big]_{j} \Big) = D_{i} A_{i}
\end{equation}
where $Z=M_{j}^{\prime}A_{i}$ and $\big[Z^{-1}A_{i}\big]_{j}$ is a supermatrix which is obtained from $Z^{-1}A_i$ by substituting $1\nu$ for the $j$-th entry.
\\
Representing entries of $A_{i}$ by $a_{i}$, one gets
\begin{equation*}
\big[Z^{-1}A_{i}\big]_{j} = \big[Z^{-1}a_{1}, \cdots, Z^{-1}a_{m}|Z^{-1}a_{m+1}, \cdots, 1\nu, \cdots, Z^{-1}a_{m+n}\big]
\end{equation*}
where $1\nu$ places in the $j_{th}$ position. Hence,
\begin{equation}\label{Z}
\big[ Z^{-1}A_{i}\big]_{j} = Z^{-1}\big[ a_{1}, \cdots, a_{m} |a_{m+1}, \cdots, Z(1\nu), \cdots, a_{m+n}\big].
\end{equation}
To have meaningful expressions, consider the following rule for $1\nu$:
\begin{equation*}
w. 1\nu = \nu (w).
\end{equation*}
Thus $Z(1\nu) = (\nu a_{j}) 1\nu=a_{j}$ and accordingly one can rewrite \eqref{Z} as below
\begin{equation*}
\big[ Z^{-1}A_{i}\big]_{j}= Z^{-1}\big[ a_{1}, \cdots, a_{m} |a_{m+1}, \cdots, a_{j}, \cdots, a_{m+n}\big] =Z^{-1}A_{i}
\end{equation*}
So for the left hand side of \eqref{g} one has:
\begin{align*}
D_{i}\Big(\big(M^{\prime}_{i}\big[Z^{-1}A_{i}\big]_{j}\big)^{-1}\big[Z^{-1}A_{i}\big]_{j}\Big) &= \Big(M^{\prime}_{i}\big(Z^{-1}A_{i}\big)\Big)^{-1} D_{i} \big(Z^{-1}A_{i}\big) = \Big( Z^{-1} \big(M^{\prime}_{i} (A_{i})\big) \Big)^{-1} Z^{-1} D_{i} A_{i}
\\
                                                                                                                                                 & = \Big( M^{\prime}_{i} (A_{i}) \Big)^{-1} Z Z^{-1} D_{i} A_{i} = D_{i}A_{i}
\end{align*}
hence $g_{ji}^{*} \circ g_{ij}^{*} = id$ holds.
\item[3)]
To prove the equality $g_{ki}^{*} \circ g_{jk}^{*} \circ g_{ij}^{*}= id$, it only suffices to show that $g_{jk}^{*} \circ g_{ij}^{*}$ is obtained from the following: 
\begin{equation*}
D_{i}\Big( \big(M^{\prime}_{i}(A_{k})\big)^{-1}A_{k} \Big) = D_{i} A_{i}
\end{equation*}
Note that the morphism $g_{jk}^{*}$ is defined by
\begin{equation}\label{M}
D_{j}\Big( \big(M^{\prime}_{j}(A_{k})\big)^{-1}A_{k} \Big) = D_{j} A_{j}
\end{equation}
To compute $g_{jk}^{*} \circ g_{ij}^{*}$, it is required to replace the supermatrix $A_{j}$ in \eqref{DJ} by the left hand side of \eqref{M}. So we have
\begin{equation} \label{final}
D_{i}\Big(\big(M^{\prime}_{i}\big[Z^{-1}A_{k}\big]_{j}\big)^{-1}\big[Z^{-1}A_{k}\big]_{j} \Big) = D_{i} A_{i}
\end{equation}
where $Z= M^{\prime}_{j}(A_{k})$.
\\
As for the second condition, one has $\big[Z^{-1}A_{k}\big]_{j}= Z^{-1}A_{k}$. Thus the left hand side of \eqref{final} is as follows:
\begin{align*}
D_{i}\Big( \big(M^{\prime}_{i}\big[Z^{-1}A_{k}\big]_{j}\big)^{-1}\big[Z^{-1}A_{k}\big]_{j} \Big) &= \Big( M^{\prime}_{i}\big(Z^{-1}A_{k}\big)\Big)^{-1}  D_{i}\big(Z^{-1}A_{k}\big) = \Big( Z^{-1} \big(M^{\prime}_{i} (A_{k})\big) \Big)^{-1}  Z^{-1} D_{i}A_{k}
\\
                                                                                                                                                    &= \Big( M^{\prime}_{i} (A_{k}) \Big)^{-1} Z Z^{-1} D_{i}A_{k} = \Big( M^{\prime}_{i} (A_{k}) \Big)^{-1} D_{i}A_{k}.
\end{align*}
By substituting in \eqref{final}, one gets
\begin{equation*}
\Big( M^{\prime}_{i} (A_{k}) \Big)^{-1} D_{i}A_{k} = D_{i}A_{i} 
\end{equation*}
So one has the following:
\begin{equation*}
 g_{jk}^{*} \circ g_{ij}^{*} = g_{ik}^{*},
\end{equation*}
This completes the proof.
\end{itemize}
\begin{flushright}
$\square$
\end{flushright}
Now, we are going to define $\widetilde{g}_{ij}$ which is introduced by \eqref{tilde}.
\begin{lemma}\label{gtilda}
Let $g^{*}_{ij}$ be the gluing morphism between sheaves of rings from $\mathcal{O}_{i}|_{U_{ij}}$ to $\mathcal{O}_{j}|_{U_{ji}}$ determined by pasting equation \eqref{DJ}. Then the following morphism is induced by $g^{*}_{ij}$ on the corresponding reduced manifolds:
\begin{equation*}
(\widetilde{g}_{ij}, \widetilde{g}^{*}_{ij}): (U_{ji}, \mathbf{C}^{\infty}(U_{j})|_{U_{ji}}) \longrightarrow (U_{ij}, \mathbf{C}^{\infty}(U_{i})|_{U_{ij}}).
\end{equation*}
\end{lemma}
\textit{Proof.}
Let $\mathcal{J}_{t}$ be the sheaf of the nilpotent elements of $\mathcal{O}_{t}$. Since $g^{*}_{ij}(\mathcal{J}_{i}) \subset \mathcal{J}_{j}$, the following morphism is induced by $g^{*}_{ij}$:
\begin{equation*}
\overline{g}^{*}_{ij}: \dfrac{\mathcal{O}_{i}}{\mathcal{J}_{i}}|_{U_{ij}} \longrightarrow \dfrac{\mathcal{O}_{j}}{\mathcal{J}_{j}}|_{U_{ji}}
\end{equation*}
On the other hand, we have the following isomorphism:
\begin{align*}
\tau_{i}&: \dfrac{\mathcal{O}_{i}}{\mathcal{J}_{i}} \longrightarrow  \mathbf{C}^{\infty}(U_{i})
\\
           &s+ \mathcal{J}_{i} \longmapsto \widetilde{s}
\end{align*}
where $\widetilde{s}(z)$, $z \in U_{i}$, is a  unique complex number so that $s-\widetilde{s}(z)$ is not invertible in $\mathcal{O}_{i}(U)$ for any neighborhood $U$ of $z$.
\\
So one has the following map between the sheaves of rings of smooth functions:
\begin{align*}
\widetilde{g}^{*}_{ij}&: \mathbf{C}^{\infty}(U_{i}) \longrightarrow  \mathbf{C}^{\infty}(U_{j})
\\
           &\widetilde{g}^{*}_{ij} := \tau_{j} \circ \overline{g}^{*}_{ij} \circ \tau^{-1}_{i}.
\end{align*}
Therefore by (\cite{Varadarajan}, Th. 4.3.1), there exists a smooth map $\widetilde{g}_{ij}: U_{j} \longrightarrow U_{i}$ such that for any $f \in \mathbf{C}^{\infty}(U_{i})$, $\widetilde{g}^{*}_{ij}(f) = f \circ \widetilde{g}_{ij}.$
\begin{flushright}
$\square$
\end{flushright}
\subsection{The underlying space of a $\nu$-projective space}
Here, we are going to identify the underlying space of a $\nu$-projective space $^{\nu}\mathcal{P}^{m|n}$.
\begin{proposition} \label{Prop}
If $(\mathcal{X}, \mathcal{O})$ be the ringed space obtained by gluing $\nu$-domains $(U_{i}, \mathcal{O}_{i})$ through $g_{ij}$, then the corresponding reduced manifold, $(\mathcal{X}, \widetilde{\mathcal{O}})$, is diffeomorphic to $\mathbb{C}P^{m}$.
\end{proposition}
To prove the theorem, we need the following lemma:
\begin{lemma}
There exists an injective immersion from $\mathcal{X}$ to $\mathbb{C}P^{m+n}$.
\end{lemma}
\textit{Proof.}
To define an injective immersion from $\mathcal{X}$ to $\mathbb{C}P^{m+n}$, it is sufficient to define a family of smooth maps $\{\psi_{it}: U_{i} \longrightarrow V_{t}\}$ where $U_{i}$ is a $\nu$-domain and $(V_{t}, \phi_{t})$, $t \in \{1, \cdots, m+n+1\}$, is the chart on $\mathbb{C}P^{m+n}$. Obviously, $V_{t}$ is in bijection with $\mathbb{C}^{m+n}$.
\\
Now, suppose that $A_{i}$ is the label of $U_{i}$, previously introduced in 3.1. Define $\psi_{it}: U_{i} \longrightarrow \mathbb{C}^{m+n}$ as a map with the component functions to be equal to the entries of the following matrix:
\begin{equation}\label{components}
D_{t} \Big( \big( M_{t} (\nu^{\prime\prime}A_{i})\big)^{-1} \nu^{\prime\prime}A_{i}\Big)
\end{equation}
where $\nu^{\prime\prime}: \mathcal{O} \longrightarrow \widetilde{\mathcal{O}} $ is a map defined by $\nu^{\prime \prime}a = \widetilde{\nu^{p(a)}a}$ on each homogeneous element $a$. In addition, for row matrix $A_{i}=(a_{ij})$ we set $\nu^{\prime\prime}A_{i} = (\nu^{\prime\prime}a_{ij})$.
\\
Let $\widetilde{g}_{ij}$ be a map introduced in lemma \ref{gtilda} and let $\theta_{ts}$ be the transition map from  $\phi_{t}$ to $\phi_{s}$. Then one has
\begin{equation*}
\theta_{ts} \circ \psi_{it}=\psi_{js} \circ \widetilde{g}_{ji}.
\end{equation*}
Note that the component functions of $\psi_{it}$ are equal to the entries of \eqref{components}. Thus $\theta_{ts} \circ \psi_{it}$ is a map whose component functions are equal to the entries of the following matrix:
\begin{equation*}
D_{s}\Big[ \Big( M_{s}\big[\big(M_{t}(\nu^{\prime\prime}A_{i})\big)^{-1}\nu^{\prime\prime}A_{i}\big]_{t}\Big)^{-1}\big[\big(M_{t} (\nu^{\prime\prime}A_{i})\big)^{-1}\nu^{\prime\prime}A_{i}\big]_{t}\Big] \\
\end{equation*}
which by the proof of the proposition \ref{gluing condition}, one gets
\begin{align*}
         &  D_{s} \Big[ \Big( M_{s}\big((M_{t}(\nu^{\prime\prime}A_{i}))^{-1}\nu^{\prime\prime}A_{i}\big)\Big)^{-1}\big(M_{t}(\nu^{\prime\prime}A_{i})\big)^{-1}\nu^{\prime\prime}A_{i}\Big]
\\
         &= D_{s}\Big[ \Big( \big(M_{t}(\nu^{\prime\prime}A_{i})\big)^{-1}\big(M_{s}(\nu^{\prime\prime}A_{i})\big)\Big)^{-1}\big(M_{t}(\nu^{\prime\prime}A_{i})\big)^{-1}\nu^{\prime\prime}A_{i}\Big] 
         \\
         &= D_{s}\Big[ \big(M_{s}(\nu^{\prime\prime}A_{i})\big)^{-1}\big(M_{t}(\nu^{\prime\prime}A_{i})\big)\big(M_{t}(\nu^{\prime\prime}A_{i})\big)^{-1}\nu^{\prime\prime}A_{i}\Big] 
         \\
         &= D_{s}\Big[ \big(M_{s}(\nu^{\prime\prime}A_{i})\big)^{-1}\nu^{\prime\prime}A_{i}\Big].
\end{align*}
On the other hand, the components of the map $\widetilde{g}_{ji}$ are equal to the entries of the matrix $\big(M_{j}(\nu^{\prime\prime}A_{i})\big)^{-1}\nu^{\prime\prime}A_{i}$. Therefore, the components of $\psi_{js} \circ \widetilde{g}_{ji}$ are the entries of the following matrix:
\begin{align*}
D_{s}&\Big[\Big(M_{s}\big((M_{j}\nu^{\prime\prime}A_{i})^{-1} \nu^{\prime\prime}A_{i}\big)\Big)^{-1} (M_{j}\nu^{\prime\prime}A_{i})^{-1} \nu^{\prime\prime}A_{i}\Big] =D_{s}\Big[\Big((M_{j}\nu^{\prime\prime}A_{i})^{-1} (M_{s}\nu^{\prime\prime}A_{i})\Big)^{-1} (M_{j}\nu^{\prime\prime}A_{i})^{-1} \nu^{\prime\prime}A_{i}\Big]
\\
         &=D_{s}\Big[(M_{s}\nu^{\prime\prime}A_{i})^{-1}(M_{j}\nu^{\prime\prime}A_{i})  (M_{j}\nu^{\prime\prime}A_{i})^{-1} \nu^{\prime\prime}A_{i}\Big] =D_{s}\Big[(M_{s}\nu^{\prime\prime}A_{i})^{-1}\nu^{\prime\prime}A_{i}\Big]
\end{align*}
Hence, $\theta_{ts} \circ \psi_{it} = \psi_{js} \circ \widetilde{g}_{ji}$.
\\
Thus the maps $\{\psi_{it}\}$ are coordinate representations of a map say $\psi$ from $\mathcal{X}$ to $\mathbb{C}P^{m+n}$.
\\
Now, we are going to show that $\psi: \mathcal{X} \longrightarrow \mathbb{C}P^{m+n}$ is an injective immersion.
\\
The smoothness of each map $\psi_{it}$ results the smoothness of $\psi$. The left inverse of $\psi$ is a map with component functions which are equal to the entries of the following matrix:
\begin{equation*}
D_{i}\Big((M_{i}\{Y\}_{t})^{-1} \{Y \}_{t}\Big)
\end{equation*}
where $\{Y\}_{t}=[y_{1}, ...y_{t-1}, 1, y_{t}, ..., y_{m}]$ whenever $Y= [y_1, ...y_{t-1}, y_t, ..., y_m].$
\\
Using $\phi$ as the left inverse of $\psi$, one has
\begin{equation*}
\phi_{ti} \circ \psi_{it} = id_{U_{i}}.
\end{equation*}
Thus, $\psi$ is an injective immersion.
\begin{flushright}
$\square$
\end{flushright}
\textit{Proof of proposition \ref{Prop}}.
Consider a map $\Lambda: \mathbb{C}P^{m} \longrightarrow \mathbb{C}P^{m+n}$ defined pointwise by $\Lambda(P)= \pi_{1}(P)$ where $\pi_{1}$ is a map induced by the following map:
\begin{align*}
\pi: &\mathbb{C}^{m} \longrightarrow \mathbb{C}^{m+n}
\\
      &(z_{1}, \cdots, z_{m}) \longmapsto (z_{1}, \cdots, z_{m}, 0, \cdots, 0)
\end{align*}
Obviously, $\Lambda$ is an imbedding and $\Lambda(\mathbb{C}P^{m}) = \psi(\mathcal{X})$. Hence, there exists a unique diffeomorphism such as $\overline{\Lambda}: \mathcal{X} \longrightarrow \mathbb{C}P^{m}$ with the property $\Lambda \circ \overline{\Lambda} = \psi$.
\begin{flushright}
$\square$
\end{flushright}
\subsection{Super line bundles}
\begin{proposition}
There exists a canonical $1|0$-super line bundle over $^{\nu}\mathcal{P}^{m|n}$.
\end{proposition}
\textit{Proof.}
On each neighborhood $U_{i}$, define a sheaf of $\mathcal{O}_{i}$-modules of rank $1|0$ as $\mathcal{O}_{i} \otimes_{\mathbb{C}} \big< A_{i} \big>_{\mathbb{C}}$, where by $\big< A_{i} \big>_{\mathbb{C}}$ we mean the super vector space generated by $A_{i}$ over $\mathbb{C}$. The sheaves may be glued together through the morphisms as below:
\begin{align*}
\eta_{ij}: &\mathcal{O}_{i}|_{U_{ij}} \otimes_{\mathbb{C}} \big<A_{i}\big>_{\mathbb{C}} \longrightarrow \mathcal{O}_{j}|_{U_{ji}} \otimes_{\mathbb{C}} \big<A_{j}\big>_{\mathbb{C}} 
\\
               & a \otimes A_{i} \longmapsto g^{*}_{ij}(a) \big(M^{\prime}_{i}(A_{j})\big)^{-1}\otimes A_{j}.
\end{align*}
It can be shown that these morphisms satisfy the conditions of the proposition \ref{gluing condition}. Therefore one obtains a $1|0$-super line bundle over $^{\nu}\mathcal{P}^{m|n}$ which we denote it by $^{\nu}\gamma_{1}$.
\begin{flushright}
$\square$
\end{flushright}
\section{Analytic approach to $\nu$-classes}
Here, we are going to generalize the exponential sheaf sequence in supergeometry, then by introducing a proper cocycle associated to $^{\nu}\mathcal{P}^{m|n}$, we define $\nu$-classes. At the end, we generalize a part of de Rham theorem and provide differential forms representing $\nu$-classes.
\subsection{The generalized exponential sheaf sequence}
Let $\mathbb{C}_{\nu_{0}}$ be the ring generated by an indeterminate $\nu_{0}$ with condition $\nu_{0}^{2}=1$. One may decompose $\mathcal{O} \otimes \mathbb{C}_{\nu_{0}}$ as follows:
\begin{equation*}
\mathcal{O} \otimes \mathbb{C}_{\nu_{0}} = (\mathcal{O}^{0} \otimes \mathbb{C}_{\nu_{0}}) \oplus (\mathcal{O}^{1} \otimes \mathbb{C}_{\nu_{0}})
\end{equation*}
where $\mathcal{O}^0$ and $\mathcal{O}^1$ are respectively the even and odd parts of the sheaf $\mathcal{O}$.
\begin{lemma}\label{exact sequence}
The following short sequence is exact.
\begin{equation*}
0 \longrightarrow \dfrac{\mathbb{Z}}{2}[\nu_{0}] \longrightarrow \mathcal{O}^{0} \oplus \nu_{0}\mathcal{O}^1 \overset{E}{\longrightarrow} \dfrac{(\mathcal{O}^{0} \otimes \mathbb{C}_{\nu_{0}})^{*}}{\{-1, +1\}} \longrightarrow 0
\end{equation*}
where $(\mathcal{O}^{0} \otimes \mathbb{C}_{\nu_{0}})^{*}$ denotes the subsheaf of even invertible elements of $(\mathcal{O}^{0} \otimes \mathbb{C}_{\nu_{0}})$. In addition the second map is defined as below
\begin{equation*}
p+ q\nu_{0} \longmapsto p+ q\nu(1)\nu_0
\end{equation*}
in which $p, q \in \dfrac{\mathbb{Z}}{2}$, and the third map may be obtained by the following map:
\begin{equation}\label{Eprime}
E^{\prime}: f + \nu_0g \longmapsto \exp \big(2\pi \sqrt{-1} (f+\nu_{0}\nu g)\big)
\end{equation}
where $f\in \mathcal{O}^{0}, g \in \mathcal{O}^{1}$  and $\exp(h)= \sum_{n=0}^{\infty}  \frac{h^{n}}{n!}$. Note that $\nu$ is an odd involution defined previously in subsection 2.2.
\end{lemma}
\textit{Proof.}
If $p+ q\nu(1)\nu_0=0$, then $p=0$ and $q=0$. Thus the second map is one-to-one.
\\
Now, we show that the kernel of the third map is equal to the image of the second map, i.e., $ker(E)=\{p+q\nu(1)\nu_{0} ;\quad  p, q \in \dfrac{\mathbb{Z}}{2}\}$. To this end, one should note that
\begin{equation*}
E^{\prime}(f + \nu_0g) = \exp \big(2\pi \sqrt{-1}(f+\nu_{0}\nu g)\big) = \big(\exp2\pi\sqrt{-1}f\big)\big(\cos(2\pi \nu g) + \nu_{0}\sqrt{-1}\sin(2\pi \nu g)\big)
\end{equation*}
where $f$ and $\nu g$ are the elements of $\mathcal{O}^{0}$.
\\
Hence, for any element in the form of $p+ q\nu(1)\nu_{0}$, in which $p,q \in \dfrac{\mathbb{Z}}{2}$, one has
\begin{equation*}
E^{\prime}(p+ q\nu(1)\nu_{0})=\exp\Big(2\pi \sqrt{-1}\big(p+q\nu_{0}\big)\Big)=\pm 1.
\end{equation*}
Conversely, let $E^{\prime}(f+\nu_{0} g) =\pm 1$. Set
\begin{equation*}
\nu g = g^{\prime}.
\end{equation*}
We have
\begin{equation*}
\exp\Big(2\pi \sqrt{-1}\big(f +\nu_{0}g^{\prime} \big)\Big) =\pm 1  \qquad  \Longrightarrow \qquad \exp \big(2\pi \sqrt{-1}f \big) \exp\big(2\pi \sqrt{-1} \nu_{0} g^{\prime}\big) =\pm 1.
\end{equation*}
By definition of the exponential map, one obtains
\begin{equation*} \label{Formula}
\exp\big( 2\pi \sqrt{-1}f \big) \Big{[} \cos(2\pi g^{\prime}) + \nu_{0} \sqrt{-1} \sin(2\pi g^{\prime}) \Big{]} =\pm 1
\end{equation*}
which concludes
\begin{equation*}
\sin(2\pi g^{\prime}) = 0 \quad \Longrightarrow \quad g^{\prime} = \frac{k_{1}}{2} (k_{1} \in \mathbb{Z}),\quad \Longrightarrow \quad g = \dfrac{k_{1}}{2}\nu(1)
\end{equation*}
Therefore,
\begin{equation*}
\exp (2\pi \sqrt{-1} f) = \pm 1 \quad \Longrightarrow \quad \cos(2\pi f) + \sqrt{-1} \sin(2\pi f) = \pm 1
\end{equation*}
which gives
\begin{equation*}
\sin(2\pi f) = 0 \qquad \Longrightarrow \qquad f = \frac{k_{2}}{2} (k_{2} \in \mathbb{Z}),
\end{equation*}
Hence, by setting $p=\dfrac{k_{2}}{2}$, $q=\dfrac{k_{1}}{2}$, one has $ker(E)= \{p+q\nu(1)\nu_{0} ;\quad  p, q \in \dfrac{\mathbb{Z}}{2}\}$ which is equal to the image of the second map of the sequence.
\\
To prove that $E$ is a surjection, let $(f+\nu_{0}g) \in \dfrac{(\mathcal{O}^{0} \otimes \mathbb{C}_{\nu_{0}})^{*}}{\{-1, +1\}}$. We should find an element of $\mathcal{O}^{0} \oplus \nu_{0}\mathcal{O}^1$, say $f^{\prime}+\nu_{0}g^{\prime}$, such that
\begin{equation*}
E^{\prime}(f^{\prime}+\nu_{0}g^{\prime}) = f+\nu_{0}g.
\end{equation*}
By definition of $E$, one gets the following equations:
\begin{equation*}
\exp(2\pi \sqrt{-1} f^{\prime})\cos(2\pi \nu g^{\prime})=f
\end{equation*}
\begin{equation*}
\exp(2\pi \sqrt{-1} f^{\prime})\sin(2\pi \nu g^{\prime})=-\sqrt{-1} g
\end{equation*}
which by considering a proper branch of logarithm, as the inverse of exponential function, one may solve these equations and find $f^{\prime}+\nu_{0}g^{\prime}$.
\begin{flushright}
$\square$
\end{flushright}
\subsection{Introduction to $\nu$-classes}
The main purpose of this subsection is to introduce a proper $1$-cocycle, $\{h_{ts}\}$, which represents an element of Cech cohomology on $\mathbb{C}P^{m}$. Then, we introduce a generalization of Chern classes in supergeometry called $\nu$-classes.
\begin{definition}
For arbitrary indices $i$ and $j$ define $h_{ij}$ as follows:
\begin{align*}
           &h_{ij} \in (\mathcal{O}^{0}_{j} \otimes \mathbb{C}_{\nu_{0}})^{*}_{|_{U_{ji}}}
\\
           &h_{ij} = \nu_{0}^{\big(p(i)+p(j)\big)}\big(M_{i}^{\prime}(A_{j})\big)^{-1}
\end{align*}
Note that $M_{i}^{\prime}(A_{j})$ is defined previously in the subsection 3.1.
\end{definition}
\begin{lemma} \label{properties of h}
$\{h_{ts}\}$ is a $1$-cocycle corresponding to the open covering $\{U_{t}\}_{1 \leq t \leq m+n+1 }$.
\end{lemma}
\textit{Proof.}
It is necessary to show that for arbitrary indices $i$, $j$ and $k$, the following equality holds on $U_{ijk}$
\begin{equation*}
 h_{jk} h_{ik}^{-1} h_{ij} = id.
\end{equation*}
Equivalently, one needs to prove the following equality on $U_{ijk}$:
\begin{equation}\label{Mprime relation}
(M_{j}^{\prime}A_{k})^{-1} M_{i}^{\prime}A_{k} = M^{\prime}_{i}A_{j}
\end{equation}
For this we consider two cases as below
\begin{itemize}
\item[1)]
If $i \leq m+1$ and $j ,k \in \{1, \cdots, m+n+1\}$, then 
\begin{equation*}
M_{i}^{\prime}A_{k} = M_{i}A_{k}, \qquad M_{i}^{\prime}A_{j} = M_{i}A_{j}.
\end{equation*}
Considering the $i$-th entries of $A_{j}$ and $A_{k}$ in the pasting equation of $g^{*}_{jk}$, i.e., $D_{j}\big( (M_{j}^{\prime}A_{k})^{-1} A_{k}\big)= D_{j}A_{j}$, one gets
\begin{equation*}
\big(M_{j}^{\prime}(A_{k})\big)^{-1} M_{i}(A_{k}) = M_{i}(A_{j}) \qquad \Longrightarrow \qquad \big(M_{j}^{\prime}(A_{k})\big)^{-1} M^{\prime}_{i}(A_{k}) = M_{i}^{\prime}(A_{j}).
\end{equation*}
\item[2)]
If $i > m+1$ and $j ,k \in \{1, \cdots, m+n+1\}$, then
\begin{equation*}
M_{i}^{\prime}(A_{k}) = \nu \big(M_{i}(A_{k})\big), \qquad M_{i}^{\prime}(A_{j}) = \nu \big(M_{i}(A_{j})\big).
\end{equation*}
In the same way as in the first case, one has $(M_{j}^{\prime}A_{k})^{-1} (M_{i}A_{k}) = M_{i}A_{j}$. So by applying $\nu$ on the two sides of this equation one gets
\begin{equation*}
(M_{j}^{\prime}A_{k})^{-1} \nu(M_{i}A_{k}) = \nu\big(M_{i}(A_{j})\big) \qquad \Longrightarrow \qquad (M_{j}^{\prime}A_{k})^{-1} M_{i}^{\prime}(A_{k}) = M_{i}^{\prime}(A_{j})
\end{equation*}
\end{itemize}
Note that for any arbitrary index $k$, $M_{k}^{\prime}(A_{k})=1$. So in the case $i=k$, one has the following:
 $$\big(M_{j}^{\prime}(A_{k})\big)^{-1} M_{k}^{\prime}(A_{k}) = M_{k}^{\prime}(A_{j}) \qquad \Longrightarrow \qquad \big(M^{\prime}_{j}(A_{k})\big)^{-1} = M^{\prime}_{k}(A_{j}).$$
Now, suppose that $i$, $j$ and $k$ be arbitrary indices. Then we have following relations on $U_{ijk}$:
\begin{align*}
& h_{jk} = \nu_{0}^{p(j) + p(k)} \big(M_{j}^{\prime}(A_{k})\big)^{-1}=\nu_{0}^{p(j) + p(k)} (M_{j}^{\prime}A_{i})^{-1} (M_{k}^{\prime}A_{i}) \in \mathcal{O}_{i_{|_{{U_{ijk}}}}}
\\
& h_{ik}= \nu_{0}^{p(k) + p(i)} \big(M_{i}^{\prime}(A_{k})\big)^{-1} = \nu_{0}^{p(k) + p(i)} \big(M_{k}^{\prime}(A_{i})\big) \in \mathcal{O}_{i_{|_{{U_{ijk}}}}}
\\
& h_{ij} = \nu_{0}^{p(i) + p(j)} \big(M_{i}^{\prime}(A_{j})\big)^{-1}=\nu_{0}^{p(i) + p(j)} \big(M_{j}^{\prime}(A_{i})\big) \in \mathcal{O}_{i_{|_{{U_{ijk}}}}}
\end{align*}
Therefore we have
\begin{align*}
h_{jk} h_{ik}^{-1}  h_{ij} &= \Big(\nu_{0}^{p(j) + p(k)} (M_{j}^{\prime}A_{i})^{-1} (M_{k}^{\prime}A_{i})\Big) \Big(\nu_{0}^{p(k) + p(i)} \big(M_{k}^{\prime}(A_{i})\big)^{-1}\Big) \Big(\nu_{0}^{p(i) + p(j)} \big(M_{j}^{\prime}(A_{i})\big)\Big)
\\
                                             &= \nu_{0}^{2\big(p(i)+p(j)+p(k)\big)} =id
\end{align*}
\begin{flushright}
$\square$
\end{flushright}
Let $\eta$ denotes the 1-cocycle $\{h_{ts}\}$. Note that $\eta$ defines an element in Cech cohomology of $\mathbb{C}P^{m}$ with coefficients in $\dfrac{(\mathcal{O}^{0} \otimes \mathbb{C}_{\nu_{0}})^{*}}{\{-1, +1\}}$ which we denote it by $\Gamma$. Hence,
\begin{equation*}
\Gamma \in H^{1}\big(\mathbb{C}P^{m}, \dfrac{(\mathcal{O}^{0} \otimes \mathbb{C}_{\nu_{0}})^{*}}{\{-1, +1\}}\big).
\end{equation*}
The short exact sequence introduced in lemma \ref{exact sequence} gives us the following long exact sequence of Cech cohomology groups:
\begin{equation*}
\cdots \longrightarrow H^{1}(\mathbb{C}P^{m}, \mathcal{O}^{0} \oplus \nu_{0}\mathcal{O}^{1}) \longrightarrow H^{1}(\mathbb{C}P^{m}, \dfrac{(\mathcal{O}^{0} \otimes \mathbb{C}_{\nu_{0}})^*}{\{-1, +1\}}) \overset{\delta}{\longrightarrow} H^{2}(\mathbb{C}P^{m},\dfrac{\mathbb{Z}}{2}[\nu_{0}])  \longrightarrow \cdots
\end{equation*}
 Define $^{\nu}c := \delta(\Gamma)$ which is called $\nu$-class of $^{\nu} \gamma_{1}$ and may be considered as an analogous of the universal Chern class in supergeometry. One may define this class for any super line bundle by a proper generalization of homotopy classiffication theorem in supergeometry.
\\
\textit{Notation:} By $H^k(M, \dfrac{\mathbb{Z}}{2}[\nu_{0}])$ we mean the $k$-th Cech cohomology group of $M$ with coefficients in $\dfrac{\mathbb{Z}}{2}[\nu_{0}]$.
\begin{theorem}
Associated to isomorphism class of each super line bundle on $\mathcal{M}=(M, \mathcal{O})$, there exists a homotopy class of a morphism $(\phi, \psi)$ from $\cal{M}$ to $^{\nu}\mathcal{P}^{m|n}$ for some pair of natural numbers, say $(m, n)$.
\end{theorem}
\textit{Proof.}
See \cite{Dr1}.
\\
Hence, one may define $\nu$-class of $\mathcal{E}$ as follows:
\begin{equation*}
^{\nu}c(\mathcal{E})=\phi^{*}(^{\nu}c) \in H^{2}\big(M, \dfrac{\mathbb{Z}}{2}[\nu_{0}]\big)
\end{equation*}
\subsection{The de Rham theorem}
Here, we are going to generalize the de Rham theorem in supergeometry for a special case. Indeed, we show that the second Cech and de Rham cohomology groups, with coefficients in $\mathbb{C}_{\nu_0}$, are isomorphic to each other. To this end, following the method used in (\cite{Griffith}, Prop.1, p.141), first, we prove lemma \ref{exact} and the theorem \ref{exact by nu_{0}}.
\\
Let $\mathcal{M}=(M, \mathcal{O})$ be a supermanifold and let $\Omega^{i}$ be a sheaf of $\mathbb{Z}_{2}$-graded super vector spaces of $i$-differential forms on the supermanifold $\mathcal{M}$ and $\Omega^{i}_{d}$ be a subsheaf of $\Omega^{i}$ whose elements are closed differential forms.
\begin{lemma}\label{exact}
The following short sequences are exact:
\begin{align}
 & 0 \longrightarrow \mathbb{R} \hookrightarrow \Omega^{0} \overset{d}{\longrightarrow} \Omega^{1}_{d} \longrightarrow 0 \label{exact seq}
 \\
 & 0 \longrightarrow \Omega^{1}_{d} \hookrightarrow \Omega^{1} \overset{d}{\longrightarrow} \Omega^{2}_{d} \longrightarrow 0 \label{exact seq2}
\end{align}
\end{lemma}
\textit{Proof.}
A sheaf sequence is exact if and only if for each $p \in M$, the corresponding sequence of stalks is exact as a sequence of super vector spaces. So for \eqref{exact seq}, it is sufficient to prove that for each $p \in M$ the sequence $0 \longrightarrow \mathbb{R}_{p} \overset{\imath}{\hookrightarrow} \Omega^{0}_{p} \overset{d}{\longrightarrow} \Omega^{1}_{d, p} \longrightarrow 0$ is exact. Equivalently, one needs to prove the following equalities:
\begin{align*}
&Ker(\imath) =0,
\\
&Im(\imath) = \mathbb{R}_{p} = Ker(d),
\\
& Im(d) = \Omega^{1}_{d, p}.
\end{align*}
The first equality is obvious.
\\
To prove the second equality, let $[f]_{p} \in \Omega^{0}_{p}$ denotes the germ of $f \in \mathcal{O}(U)$ such that $d[f]_{p} = [0] \in \Omega^{1}_{d, p}$. Then, there exists a small enough contractible neighborhood $V \subset U$ of $p$ such that it splits on $M$ and one has $df_{|_{V}} = 0$. Therefore, by (\cite{Kostant}, Th. 4.6), $f$ is constant, i.e., $[f]_{p} \in \mathbb{R}_{p}$. Conversely, let $[f]_{p} \in \mathbb{R}_{p}$ then $d[f]_{p} =0$. So $Im(\imath) = \mathbb{R}_{p} = Ker(d)$.
\\
By (\cite{Kostant}, Th. 4.6), $d: \Omega^{0}_{p} \longrightarrow \Omega^{1}_{d, p}$ is a surjective map, so $Im(d) = \Omega^{1}_{d, p}$.
\\
Similarly, one may prove that \eqref{exact seq2} is an exact sequence.
\begin{flushright}
$\square$
\end{flushright}
\begin{theorem}\label{exact by nu_{0}}
The following short sequences are exact:
\begin{align}
 & 0 \longrightarrow \mathbb{R} \otimes \mathbb{C}_{\nu_{0}} \hookrightarrow \Omega^{0} \otimes \mathbb{C}_{\nu_{0}} \overset{d}{\longrightarrow} \Omega^{1}_{d} \otimes \mathbb{C}_{\nu_{0}} \longrightarrow 0 \label{exact seq tensor}
 \\
 & 0 \longrightarrow \Omega^{1}_{d} \otimes \mathbb{C}_{\nu_{0}} \hookrightarrow \Omega^{1} \otimes \mathbb{C}_{\nu_{0}} \overset{d}{\longrightarrow} \Omega^{2}_{d} \otimes \mathbb{C}_{\nu_{0}} \longrightarrow 0 \label{exact seq2 tensor}
\end{align}
\end{theorem}
\textit{Proof.}
The exactness of the sequence \eqref{exact seq} at the right shows that the short sequence \eqref{exact seq tensor} is exact at the right. In addition, $\mathbb{R} \otimes \mathbb{C}_{\nu_{0}} \hookrightarrow \Omega^{0} \otimes \mathbb{C}_{\nu_{0}}$ is an injection.
\\
Similarly, one may prove the exactness of the sequence \eqref{exact seq2 tensor}.
\begin{flushright}
$\square$
\end{flushright}
This theorem leads us to generalize the de Rham theorem as below.
\\
For an arbitrary sheaf $\mathcal{O}$ on a supermanifold $M$ and a locally finite open cover $\underline{U}=\{U_{\alpha}\}$, let $C^{p}(\underline{U}, \mathcal{O}) = \underset{\alpha_{0} \neq \alpha_{1} \neq \cdots \neq \alpha_{p}}{\prod} \mathcal{O}(U_{\alpha_{0}} \cap \cdots \cap U_{\alpha_{p}})$ denotes the set of $p$-cochains of $\mathcal{O}$. Then we have the following short exact sequence of cochain groups:
\begin{equation*}
0 \longrightarrow C^{p}(\underline{U}, \mathbb{C}_{\nu_{0}}) \overset{\imath}{\longrightarrow} C^{p}(\underline{U}, \Omega^{0} \otimes \mathbb{C}_{\nu_{0}}) \overset{d}{\longrightarrow} C^{p}(\underline{U}, \Omega^{1}_{d} \otimes \mathbb{C}_{\nu_{0}}) \longrightarrow 0
\end{equation*}
which by \cite{Hirzebruch}, it gives rise to an associated long exact sequence of cohomology groups:
\begin{align}\label{Cohomology sequence}
 0 &\longrightarrow H^{0}(M, \mathbb{C}_{\nu_{0}}) \longrightarrow H^{0}(M, \Omega^{0}\otimes \mathbb{C}_{\nu_{0}}) \longrightarrow H^{0}(M, \Omega^{1}_{d}\otimes \mathbb{C}_{\nu_{0}}) \nonumber \\ 
    &\longrightarrow H^{1}(M, \mathbb{C}_{\nu_{0}}) \longrightarrow H^{1}(M, \Omega^{0}\otimes \mathbb{C}_{\nu_{0}}) \longrightarrow H^{1}(M, \Omega^{1}_{d}\otimes \mathbb{C}_{\nu_{0}}) \\
    & \longrightarrow H^{2}(M, \mathbb{C}_{\nu_{0}}) \longrightarrow H^{2}(M, \Omega^{0}\otimes \mathbb{C}_{\nu_{0}}) \longrightarrow H^{2}(M, \Omega^{1}_{d}\otimes \mathbb{C}_{\nu_{0}}) \longrightarrow \cdots \nonumber
\end{align}
Note that we replace $\mathbb{R}\otimes \mathbb{C}_{\nu_{0}}$ by $\mathbb{C}_{\nu_{0}}$ because of the equality  $\mathbb{R}\otimes \mathbb{C}_{\nu_{0}}=\mathbb{C}_{\nu_{0}}$.
\\
In the same way, one has a long exact sequence of cohomology groups associated to the short exact sequence \eqref{exact seq2 tensor}.
\begin{proposition} \label{De Rham}
The following map is an isomorphism:
\begin{equation*}
\dfrac{H^{0}(M, \Omega^{2}_{d}\otimes \mathbb{C}_{\nu_{0}})}{dH^{0}(M, \Omega^{1}\otimes \mathbb{C}_{\nu_{0}})} \longrightarrow H^{2}(M, \mathbb{C}_{\nu_{0}})
\end{equation*}
\end{proposition}
\textit{Proof.}
Since $\Omega^{0}\otimes \mathbb{C}_{\nu_{0}}$ is a locally free sheaf, by (\cite{Manin}, page 188) it is fine and one has
\begin{equation*}
H^{1}(M, \Omega^{0}\otimes \mathbb{C}_{\nu_{0}}) = 0, \qquad H^{2}(M, \Omega^{0}\otimes \mathbb{C}_{\nu_{0}})=0.
\end{equation*}
Hence, a part of the sequence \eqref{Cohomology sequence} reduces to the following exact sequence:
\begin{equation*}
0 \longrightarrow H^{1}(M, \Omega^{1}_{d}\otimes \mathbb{C}_{\nu_{0}}) \overset{\delta_{1}}{\longrightarrow} H^{2}(M, \mathbb{C}_{\nu_{0}}) \longrightarrow 0
\end{equation*}
thus $\delta_{1}$ is an isomorphism.
\\
Similarly, we have the following long exact sequence associated to the short exact sequence \eqref{exact seq2 tensor}:
\begin{equation*}
0 \longrightarrow H^{0}(M, \Omega^{1}_{d}\otimes \mathbb{C}_{\nu_{0}}) \longrightarrow H^{0}(M, \Omega^{1}\otimes \mathbb{C}_{\nu_{0}}) \overset{d}{\longrightarrow} H^{0}(M, \Omega^{2}_{d}\otimes \mathbb{C}_{\nu_{0}}) \overset{\delta_{2}}{\longrightarrow} H^{1}(M, \Omega^{1}_{d}\otimes \mathbb{C}_{\nu_{0}}) \longrightarrow 0
\end{equation*}
Because of the exactness, we have $ker (\delta_{2}) = Im(d)$ which results that 
$\overline{\delta}_{2}:  \dfrac{H^{0}(M, \Omega^{2}_{d}\otimes \mathbb{C}_{\nu_{0}})}{dH^{0}(M, \Omega^{1}\otimes \mathbb{C}_{\nu_{0}})} \longrightarrow H^{1}(M, \Omega^{1}_{d}\otimes \mathbb{C}_{\nu_{0}})$
 is an isomorphism.
\\
So $\delta_{1}\overline{\delta}_{2}$ is an isomorphism and one has
\begin{equation*}
\dfrac{H^{0}(M, \Omega^{2}_{d}\otimes \mathbb{C}_{\nu_{0}})}{dH^{0}(M, \Omega^{1}\otimes \mathbb{C}_{\nu_{0}})} \cong H^{2}(M, \mathbb{C}_{\nu_{0}})
\end{equation*}
\begin{flushright}
$\square$
\end{flushright}
\begin{corollary}
Up to isomorphism, we may consider $^{\nu} c$ as an element of the de Rham cohomology.
\end{corollary}
\subsection{Computing $\nu$-classes}
Here, we are going to obtain a preimage of the cocycle $\eta=\{h_{ij}\}$ under the map $E^{\prime}$ defined in lemma \ref{exact sequence}. In fact, we should find $(f_{ij}+\nu_{0} g_{ij}) \in \mathcal{O}^{0} \oplus \nu_{0}\mathcal{O}^1$ such that 
\begin{equation}\label{preimagehij}
E^{\prime}(f_{ij}+\nu_{0} g_{ij})= h_{ij}.
\end{equation}
To this end, consider two cases as below:
\begin{itemize}
\item[1)]
If $1 \leq i, j \leq m+1$ or $i, j > m+1$, then $\nu_{0}^{p(i)+p(j)} = 1$, so we have
\begin{equation*}
h_{ij} = \big( M^{\prime}_{i}(A_{j})\big)^{-1}.
\end{equation*}
By substituting in \eqref{preimagehij}, one has
\begin{equation*}
E^{\prime}(f_{ij}+\nu_{0} g_{ij}) = \big( M^{\prime}_{i}(A_{j})\big)^{-1}
\end{equation*}
which gives
\begin{equation*}
\exp\big((2\pi \sqrt{-1})f_{ij}\big)\big(\cos(2\pi \nu g_{ij}) + \nu_{0}\sqrt{-1}\sin(2\pi \nu g_{ij})\big) = \big( M^{\prime}_{i}(A_{j})\big)^{-1}.
\end{equation*}
Since there is no coefficient of $\nu_{0}$ on the right hand side, one concludes
\begin{equation*}
\sin(2\pi \nu g_{ij}) = 0 \qquad \Longrightarrow \qquad \nu g_{ij} = 0, \qquad \cos (2\pi \nu g_{ij}) =1.
\end{equation*}
So one has
\begin{equation*}
                            \exp\big((2\pi \sqrt{-1})f_{ij}\big) = \big( M^{\prime}_{i}(A_{j})\big)^{-1}  \qquad \Longrightarrow \qquad f_{ij}= \dfrac{1}{2\pi \sqrt{-1}} \log \Big(\big( M^{\prime}_{i}(A_{j})\big)^{-1}\Big).
\end{equation*}
Hence,
\begin{equation}\label{f+nu g standard}
f_{ij}+\nu_{0}g_{ij} = \dfrac{1}{2\pi \sqrt{-1}} \log \Big(\big( M^{\prime}_{i}(A_{j})\big)^{-1}\Big).
\end{equation}
\item[2)]
If $1 \leq i \leq m+1$ and $j > m+1$, then $\nu_{0}^{p(i)+p(j)} = \nu_{0}$, so we have
\begin{equation*}
h_{ij} = \nu_{0} \big( M^{\prime}_{i}(A_{j})\big)^{-1}.
\end{equation*}
By substituting in \eqref{preimagehij}, one gets
\begin{equation}\label{preimage}
\exp\big((2\pi \sqrt{-1})f_{ij}\big)\big(\cos(2\pi \nu g_{ij}) + \nu_{0}\sqrt{-1}\sin(2\pi \nu g_{ij})\big) = \nu_{0} \big( M^{\prime}_{i}(A_{j})\big)^{-1}
\end{equation}
which yields the following relations:
\begin{equation*}
\cos(2\pi \nu g_{ij}) = 0 \qquad \Longrightarrow \qquad 2\pi \nu g_{ij}= \dfrac{\pi}{2} \qquad \Longrightarrow \qquad g_{ij}=\dfrac{1}{4} \nu(1), \qquad \sin(2\pi \nu g_{ij}) = 1.
\end{equation*}
So one may rewrite \eqref{preimage} as follows:
\begin{equation*}
\exp\big((2\pi \sqrt{-1})f_{ij}\big)\big( \nu_{0}\sqrt{-1}\big) = \nu_{0} \big( M^{\prime}_{i}(A_{j})\big)^{-1}
\end{equation*}
which concludes
\begin{align*}
                         &\exp\big((2\pi \sqrt{-1})f_{ij}\big) = -\sqrt{-1} \big( M^{\prime}_{i}(A_{j})\big)^{-1}
                         \\ \Longrightarrow \qquad
                         &f_{ij} = \dfrac{1}{2\pi \sqrt{-1}} \log\Big(-\sqrt{-1} \big( M^{\prime}_{i}(A_{j})\big)^{-1}\Big).
\end{align*}
Hence,
\begin{equation}\label{f+nu g nonstandard}
f_{ij}+\nu_{0} g_{ij} = \dfrac{1}{2\pi \sqrt{-1}} \log\Big(-\sqrt{-1} \big( M^{\prime}_{i}(A_{j})\big)^{-1}\Big) + \nu_{0}\dfrac{1}{4} \nu(1).
\end{equation}
\end{itemize}
Thus a right inverse of $E^{\prime}$, say $L$, is as follows:
\begin{equation*}
L(h_{ij}) = f_{ij} + \nu_{0} g_{ij},
\end{equation*}
where $f_{ij}$ and $g_{ij}$ are introduced by one of the equations \eqref{f+nu g standard} and \eqref{f+nu g nonstandard}.
\\
Since $\widetilde{M^{\prime}_{i}(A_{j})} \in \mathbb{C} \setminus \{0\}$, one may decompose $\nu$-domains $U_{j}$ by a decomposition of $\mathbb{C} \setminus \{0\}$.
\\
Let $V^{1}$(resp. $V^{2}$) be an open subset of $\mathbb{C}$ obtained by removing $0$ and all positive (resp. negative) real numbers from the complex plane. So one may obtain the following decomposition of $U_{j}$:
\begin{equation*}
U_{j}= U_{j}^{1} \cup U_{j}^{2}
\end{equation*}
where $U_{j}^{1}$ and $U_{j}^{2}$ are respectively the preimages of $V^{1}$ and $V^{2}$ under $h_{ij}$.
\\
Now, one may compute $L(h_{ij})$ by considering the restriction $h_{ij}|_{U_{i}^{k} \cap U_{j}^{l}}$  and a proper branch of logarithm on $U_{i}^{k} \cap U_{j}^{l}$ where $k, l=1,2$.
\\
Now, if we set
\begin{equation*}
(\delta \eta)_{ijk} =  L (h_{jk}) - L (h_{ik}) + L (h_{ij}),
\end{equation*}
then $\{(\delta \eta)_{ijk}\}$ is a representative of $^{\nu}c$.
\begin{example} \label{example of super}
Consider the super line bundle over $^{\nu} \mathcal{P}^{2|1}$. From the subsection 3.1, we know that $^{\nu}\mathcal{P}^{2|1}$ is constructed by gluing $\nu$-domains $(U_{i}, \mathcal{O}_{i}) = (\mathbb{C}^{m}, \mathcal{C}_{\mathbb{C}^{m}})$, $(1 \leq i \leq 4)$. Let $U_{i}$ be labeled by $A_{i}$ where
\begin{equation*}
\begin{array}{rc}
A_{1}=(1, z_{1}^{(1)}, z_{2}^{(1)} | e_{1}^{(1)}),
\end{array}
\\
\begin{array}{rc}
A_{2}=(z_{1}^{(2)}, 1, z_{2}^{(2)} | e_{1}^{(2)}),
\end{array}
\end{equation*}
\begin{equation*}
\begin{array}{rc}
A_{3}=(z_{1}^{(3)}, z_{2}^{(3)}, 1 | e_{1}^{(3)}),
\end{array}
\\
\begin{array}{rc}
\quad A_{4}=(z_{1}^{(4)}, z_{2}^{(4)}, \nu(e_{1}^{(4)}) | 1\nu).
\end{array}
\end{equation*}
Note that $U_{1}$, $U_{2}$  and $U_{3}$ are standard $\nu$-domains and $U_{4}$ is a nonstandard $\nu$-domain.
\\
By definition of $h_{ij}$, we have
\begin{equation*}
h_{21}=\dfrac{1}{z_{1}^{(1)}}
\qquad
h_{32}=\dfrac{1}{z_{2}^{(2)}}
\qquad
h_{43}=\nu_{0} \dfrac{1}{\nu e_{1}^{(3)}}
\qquad
h_{14}=\nu_{0} \dfrac{1}{z_{1}^{(4)}}.
\end{equation*}
Now, we are going to obtain $L(h_{21})$ under the map $E^{\prime}$ defined by \eqref{Eprime}. Since $U_{1}$ and $U_{2}$ are standard, we should apply the equation \eqref{f+nu g standard} to compute $L(h_{21})$. Hence, we have
\begin{equation*}\label{L(h)}
L(h_{21}) = \dfrac{1}{2\pi \sqrt{-1}} \log \Big(\big( M^{\prime}_{2}(A_{1})\big)^{-1}\Big) = \dfrac{1}{2\pi \sqrt{-1}} \log \big(\dfrac{1}{z_{1}^{(1)}} \big)
\end{equation*}
In this case, we may have the decomposition $U_{1}=U_{1}^{1} \cup U_{1}^{2}$ where $U_{1}^{1}$ consists of those points of $U_{1}$ satisfying in at least one of the following relations:
\begin{itemize}
\item[1)]
$real(z_{1}^{(1)}) < 0$,
\item[2)]
$Im(z_{1}^{(1)}) \neq 0$,
\end{itemize}
and similarly $U_{1}^{2}$ consists of those points of $U_{1}$ for which at least one of the relations $real(z_{1}^{(1)}) > 0$ and $Im(z_{1}^{(1)}) \neq 0$ holds.
\\
In the same way, one may obtain such a decomposition for the $\nu$-domain $U_{2}$.
\\
Considering the restriction of $h_{21}|_{{U_{1}^{1} \cap U_{2}^{1}}}$, one has $h_{21} = \dfrac{1}{z_{1}^{(1)}}$ such that $real(z_{1}^{(1)}) < 0$ or $Im(z_{1}^{(1)}) \neq 0$.
\\
So by assuming $0<arg(z_{1}^{(1)})<2\pi$, one has
\begin{align*}
L(h_{21}) &= \dfrac{1}{2\pi \sqrt{-1}} \log \big(\dfrac{1}{z_{1}^{(1)}} \big) = \dfrac{1}{2\pi \sqrt{-1}} \log \big(\dfrac{\overline{z_{1}^{(1)}}}{\vert z_{1}^{(1)}\vert ^{2}} \big)
\\
                                                  &= \dfrac{1}{2\pi \sqrt{-1}} \big( -\log \vert z_{1}^{(1)}\vert + \sqrt{-1}(2\pi - arg(z_{1}^{(1)})\big).
\end{align*}
By similar computations, one may obtain a representative of $^{\nu}c$ as follows:
\begin{equation*}
(\delta \eta)_{241} = L(h_{41}) - L(h_{21}) + L(h_{24}) = \dfrac{-1}{2}
\end{equation*}
\end{example}
\subsection{De Rham representative  of the $\nu$-class for super line bundles}
In this subsection, we are going to find a closed 2-form representing $\nu$-classes. In \cite{Griffith}, it is proved that for a complex line bundle $\gamma$ on a compact manifold $M$ of dimension $n$, there exists a representative of Chern classes in terms of the curvature form $\Theta$ as follows:
\begin{equation*}
c = \Big[\dfrac{1}{2\pi \sqrt{-1}}\Theta\Big] \in H^{2}_{DR}(M).
\end{equation*}
Now, we introduce a curvature form $\Theta$ in terms of $\{h_{ij}\}$ by applying a partition of unity. Inspired by it, we obtain the desired closed 2-form representing $\nu$-classes.
\\
Let $M$ be a compact manifold of dimension $n$ and $\gamma$ denotes a complex line bundle on it. Let $\{\rho_{\alpha}\}$ be a partition of unity subordinate to $\{U_{\alpha}\}$, an open cover of $M$. A connection 1-form $\theta_{\alpha}$ may be defined as follows:
\begin{align*}
\theta_{\alpha} s_{\alpha} &= \nabla (s_{\alpha}) = \sum_{\beta} \rho_{\beta} \nabla^{\beta} (s_{\alpha}) 
\\
                                       &= \sum_{\beta} \rho_{\beta} \nabla^{\beta} (h_{\alpha \beta} s_{\beta}) = \sum_{\beta} \rho_{\beta} d(h_{\alpha \beta}) s_{\beta},
\end{align*}
where $s_{\alpha}$ is a section on $U_{\alpha}$ such that $s_{\alpha}(x)$ is a basis for the fiber of $\gamma$ over $x$. By restriction to $U_{\alpha} \cap U_{\alpha^{\prime}}$, one gets
\begin{equation*}
\theta_{\alpha}h_{\alpha \alpha^{\prime}} s_{\alpha^{\prime}} = \sum_{\beta} \rho_{\beta} d(h_{\alpha \beta}) h_{\beta \alpha^{\prime}} s_{\alpha^{\prime}},
\end{equation*}
hence,
\begin{equation}\label{theta_{i}}
\theta_{\alpha} = \sum_{\beta} \rho_{\beta} d\big(\log(h_{\alpha \beta})\big).
\end{equation}
So we may obtain the curvature form $\Theta$ as follows:
\begin{equation*}\label{dtheta}
\Theta = d\theta_{\alpha} = \sum_{\beta} d(\rho_{\beta}) d\big(\log(h_{\alpha \beta})\big).
\end{equation*}
Note that the cocycle $\eta$ defined in the subsection 4.2, does not correspond to any super line bundle. So by applying formula \eqref{theta_{i}} for the cocycle $\eta$ and substituting $L(h_{ij})$ instead of $\log (h_{ij})$, we may obtain local 1-forms which their exterior derivatives define a closed 2-form representing $^{\nu}c$.
\begin{theorem}
For a super line bundle $^{\nu}\gamma_{1}$, there exists a representation of $\nu$-class $^{\nu}c$ by a closed 2-form corresponding to the cocycle $\eta$.
\end{theorem}
\textit{Proof.}
Inspired by the previous part, we may define a 1-form as $\omega_{i} = \sum_{j} \rho_{j} d\big(L (h_{ij})\big)$ on $U_{i}$ which its exterior derivative is the following 2-form:
\begin{equation*}\label{Omega}
R_{i} = \sum_{j} d(\rho_{j}) d\big(L (h_{ij})\big) =\dfrac{1}{2\pi \sqrt{-1}} \sum_{j}\dfrac{1}{h_{ij}} d(\rho_{j}) d(h_{ij}) + \nu_{0} \sum_{j} \dfrac{k_{ij}}{4} d(\rho_{j}) d\big(\nu(1)\big)
\end{equation*}
where $k_{ij}=0$ if $1 \leq i, j \leq m+1$ or $i, j > m+1$, otherwise $k_{ij}=1$.
\\
Now, we prove that $R_{i}$ defines a global closed 2-form $R$. To this end, we compute $\omega_{i} - \omega_{j}$.
\\
Let $1 \leq i \leq m+1$ and $j > m+1$. So one has
\begin{align*}
\omega_{i} - \omega_{j} &= \sum_{k=1}^{m+n+1} \rho_{k} d\big(L (h_{ik})\big) - \sum_{k=1}^{m+n+1} \rho_{k} d\big(L (h_{jk})\big)
\\
                                      &= \sum_{k=1}^{m+1} \rho_{k} d\Big( \dfrac{1}{2\pi \sqrt{-1}} \log\big(M_{i}^{\prime}(A_{k})\big)^{-1} \Big) + \sum_{k=m+2}^{m+n+1} \rho_{k} d\Big( \dfrac{1}{2\pi \sqrt{-1}} \log\big(-\sqrt{-1} \big(M_{i}^{\prime}(A_{k})\big)^{-1} \big) +\nu_{0} \dfrac{1}{4} \nu(1)\Big) \\
                                      &-\sum_{k=1}^{m+1} \rho_{k} d\Big( \dfrac{1}{2\pi \sqrt{-1}} \log\big(-\sqrt{-1} \big(M_{j}^{\prime}(A_{k})\big)^{-1} \big) +\nu_{0} \dfrac{1}{4} \nu(1)\Big) - \sum_{k=m+2}^{m+n+1} \rho_{k} d\Big( \dfrac{1}{2\pi \sqrt{-1}} \log\big(M_{j}^{\prime}(A_{k})\big)^{-1} \Big) 
                                      \\
                                      &= \dfrac{1}{2\pi \sqrt{-1}} \sum_{k=1}^{m+1} \rho_{k} d\Big(\log\big(\sqrt{-1} \big(M_{i}^{\prime}(A_{k})\big)^{-1} \big(M_{j}^{\prime}(A_{k})\big) \big) - \nu_{0} \dfrac{m+1}{4} d\big(\nu(1)\big) 
                                      \\
                                      &+ \dfrac{1}{2\pi \sqrt{-1}} \sum_{k=m+2}^{m+n+1} \rho_{k} d\Big( \log\big(-\sqrt{-1}\big(M_{i}^{\prime}(A_{k})\big)^{-1} \big(M_{j}^{\prime}(A_{k})\big)\big)\Big) + \nu_{0} \dfrac{n}{4} d\big(\nu(1)\big)
\end{align*}
hence, by the equality \eqref{Mprime relation} proved in lemma \ref{properties of h}, we have
\begin{align*}
\omega_{i} - \omega_{j}
                                     &= \dfrac{1}{2\pi \sqrt{-1}} d\Big( \log\big(-\sqrt{-1}\big(M_{i}^{\prime}(A_{j})\big)^{-1}\big)\Big) + \nu_{0} \dfrac{n-(m+1)}{4} d\big(\nu(1)\big) 
                                     \\
                                     &= d\big(L(h_{ij})\big) + \nu_{0} \dfrac{n-(m+1)-1}{4} d\big(\nu(1)\big)
\end{align*}
Therefore, we have
\begin{equation*}
\omega_{i} - \omega_{j} = d \big(L(h_{ij})\big) + \nu_{0} \dfrac{K_{ij}-1}{4} d\big(\nu(1)\big)
\end{equation*}
where $K_{ij}=n-(m+1)$ if $1 \leq i \leq m+1$ and $m+1 < j$ and $K_{ij} = 1$ if $1 \leq i,j \leq m+1$ or $m+2 \leq i, j \leq m+n+1$.
\\
From the relation $\omega_{i} - \omega_{j} = d\Big(L(h_{ij}) + \nu_{0} \dfrac{K_{ij}-1}{4} \nu(1)\Big)$, we deduce that $(d\omega_{i})_{|_{U_{ij}}} = (d\omega_{j})_{|_{U_{ij}}}$, that is $(R_{i})_{|_{U_{ij}}}=(R_{j})_{|_{U_{ij}}}$ which by the globality property of sheaves it proves that we have a global closed 2-form $R$ defined by $R_{|_{U_{i}}} := R_{i}$.
\\ 
Hence, by isomorphisms $\delta_{1}$ and $\overline{\delta}_{2}$ defined in proposition \ref{De Rham}, one has
\begin{align*}
&\big(\overline{\delta}_{2}(R)\big)_{j, i}= dL (h_{ij}) + \nu_{0} \dfrac{K_{ij}-1}{4} d\big(\nu(1)\big),
\\
&\big(\delta_{1}\overline{\delta}_{2}(R)\big)_{k, l, r} =\Big(\delta_{1} d \big(L (h_{ij}) + \nu_{0} \dfrac{K_{ij}-1}{4} \nu(1)\big)\Big)_{k, l, r} 
\\
 & \qquad  \qquad  \qquad = L(h_{lr}) -L (h_{kr}) + L (h_{kl}) + \nu_{0} \dfrac{(K_{lr}-1)-(K_{kr}-1)+(K_{kl} -1)}{4} \nu(1)
\end{align*}
A straightforward computation shows that $\dfrac{(K_{lr}-1)-(K_{kr}-1)+(K_{kl} -1)}{4}=0$ for any arbitrary indices $k, l, r$. Hence, 
\begin{equation*}
\big(\delta_{1}\overline{\delta}_{2}(R)\big)_{k, l, r} =\Big(\delta_{1} d \big(L (h_{ij}) + \nu_{0} \dfrac{K_{ij}-1}{4} \nu(1)\big)\Big)_{k, l, r} = L(h_{lr}) -L (h_{kr}) + L (h_{kl})
\end{equation*}
Therefore, one may represent $^{\nu}c$ in terms of the closed 2-form $R$.
\begin{flushright}
$\square$
\end{flushright}
\subsection{Representation of $\nu$-classes for super vector bundles}
In the previous section, we show that the (first) $\nu$-class of the canonical $1|0$-super line bundle over $\nu$-projective space may be represented by a closed 2-form in $\Omega\otimes \mathbb{C}_{\nu_{0}}$. So one may deduce that $\nu$-classes of higher rank may be represented by higher degree closed forms. In this section, in a similar way of representing universal Chern classes in common geometry, for the canonical super vector bundle over $\nu$-grassmannian $^{\nu}Gr_{k|l}(m|n)$ \cite{Bahadory}, we introduce closed forms in $\Omega_{\mathcal{G}}\otimes \mathbb{C}_{\nu_{0}}$ which may be considered as representatives of $\nu$-classes of rank higher than 1. By $\Omega_{\mathcal{G}}$ we mean the structure sheaf of $\nu$-grassmannian. In \cite{Husemoller}, Chern classes of a complex vector bundle of rank $k$ are represented by closed $2r$-forms($1 \leq r \leq k$) which are polynomials in the curvature form of a connection on the bundle. One may obtain the curvature 2-form matrix $\Theta$ for a certain connection by a method as in the previous subsection.
\\
Let $M$ denotes a manifold of dimension $n$ and $\xi^{k}$ a complex vector bundle over it. There exists a partition of unity $\{\rho_{\alpha}\}$ subordinate to an open covering $\{U_{\alpha}\}$ of $M$. Let $\{h^{\alpha\beta}\}$ be coordinate transformations of $\xi^{k}$ and $\nabla$ be a connection defined as follows: 
\begin{align*}
\nabla(s_{i}^{\alpha}) &= \sum_{\beta}\rho_{\beta} \nabla^{\beta}(s_{i}^{\alpha})= \sum_{\beta}\rho_{\beta} \nabla^{\beta}(\sum_{l} h_{il}^{\alpha \beta} s_{l}^{\beta})
\\
                                  &= \sum_{l}\sum_{\beta} \rho_{\beta} d(h_{il}^{\alpha \beta}) s_{l}^{\beta} = \sum_{j}\sum_{l}\sum_{\beta} \rho_{\beta} d(h_{il}^{\alpha \beta}) h_{lj}^{\beta\alpha }s_{j}^{\alpha},
\end{align*}
where $s_{i}^{\alpha}$, ($1\leq i \leq k$), are sections of $\xi|_{U_{\alpha}}$ such that for each $x\in U_{\alpha}$, $\{s_{i}^{\alpha}(x)\}$ is a basis for the fiber over $x$.
\\
Let $\theta^\alpha= (\theta^{\alpha}_{ij})$ be the connection form of $\nabla$. Then, $\nabla(s_{i}^{\alpha}) = \sum_{j} \theta_{ij}^{\alpha}s_{j}^{\alpha}$.
\\
Therefore, we have
\begin{equation*}
\theta_{ij}^{\alpha} = \sum_{l}\sum_{\beta} \rho_{\beta} d(h_{il}^{\alpha \beta}) h_{lj}^{\beta\alpha},
\end{equation*}
or one may write
\begin{equation*}\label{connection form}
\theta^{\alpha} = \sum_{\beta} \rho_{\beta} d(h^{\alpha \beta}) h^{\beta\alpha}
\end{equation*}
By a straightforward calculation, one may prove that $\Theta^{\alpha} = d\theta^{\alpha} - \theta^{\alpha} \wedge \theta^{\alpha}$ is the corresponding curvature 2-form matrix. So we have
\begin{align*}
\Theta^{\alpha} &= d\theta^{\alpha} - \theta^{\alpha} \wedge \theta^{\alpha} 
\\
                         &= \sum_{\beta} d(\rho_{\beta}) d(h^{\alpha \beta}) h^{\beta\alpha} - \sum_{\beta} \rho_{\beta} d(h^{\alpha \beta}) d(h^{\beta\alpha}) - \sum_{\beta, \beta^{\prime}}\rho_{\beta} d(h^{\alpha \beta}) h^{\beta\alpha} \wedge \rho_{\beta^{\prime}} d(h^{\alpha \beta^{\prime}}) h^{\beta^{\prime}\alpha}
\end{align*}
Now, we apply this technique to obtain a supermatrix-valued 2-form for a canonical $k|l$-super vector bundle $^{\nu} \gamma_{k|l}$ over $\nu$-grassmannian $^{\nu}Gr_{k|l}(m|n)$. To obtain such a 2-form, we should define a cocycle $\{h^{\alpha\beta}\}$. By considering gluing morphisms of $^{\nu} \gamma_{k|l}$ introduced in \cite{Bahadory}, one may define
\begin{equation}\label{cocycle of grassmannian}
h^{\alpha\beta} = \nu_{0}^{p(\alpha)+p(\beta)} (M_{\alpha}^{\prime}A_{\beta})^{-1}
\end{equation}
where an index $\alpha=I|J$ is a multi-index such that $I$ and $J$ are respectively ordered sets of $\{1, \cdots, m\}$ and $\{1, \cdots, n\}$ with the property that $|I|+|J|=k+l$. In addition, $p(\alpha)=0$ if $|I|=k$, otherwise $p(\alpha)=1$. Moreover, $A_\alpha$ is a $k|l \times m|n$ supermatrix associated to the $\nu$-domain $U_\alpha$. By $M^{\prime}_{\alpha}A_{\beta}$ we mean a $k|l \times k|l$ standard supermatrix which is obtained by applying certain modifications on columns of $A_\beta$ with indices in $\alpha$. For more information, see (\cite{Bahadory}, section 2).
\\
One can easily show that $\{h^{\alpha\beta}\}$ defined by \eqref{cocycle of grassmannian} is a cocycle. Now, consider a partition of unity $\{\rho_{\alpha}\}$ subordinate to the open covering $\{U_{\alpha}\}$ of $^{\nu}Gr_{k|l}(m|n)$. We define a ($k|l \times k|l$) supermatrix-valued 1-form as follows:
\begin{equation} \label{omega def}
\omega^{\alpha} = \sum_{\beta} \rho_{\beta} d(h^{\alpha \beta}) h^{\beta\alpha}
\end{equation}
By setting $R^{\alpha} = d\omega^{\alpha} - \omega^{\alpha} \wedge \omega^{\alpha}$, we obtain a ($k|l \times k|l$) supermatrix-valued 2-form as follows:
\begin{equation} \label{Omega defi}
R^{\alpha}= \sum_{\beta} d(\rho_{\beta}) d(h^{\alpha \beta}) h^{\beta\alpha} - \sum_{\beta} \rho_{\beta} d(h^{\alpha \beta}) d(h^{\beta\alpha}) - \sum_{\beta, \beta^{\prime}}\rho_{\beta} d(h^{\alpha \beta}) h^{\beta\alpha} \wedge \rho_{\beta^{\prime}} d(h^{\alpha \beta^{\prime}}) h^{\beta^{\prime}\alpha}.
\end{equation}
One may prove that $\omega^{\alpha^{\prime}} = d(h^{\alpha^{\prime}\alpha}) (h^{\alpha^{\prime}\alpha})^{-1} + (h^{\alpha^{\prime}\alpha}) \omega^{\alpha} (h^{\alpha^{\prime}\alpha})^{-1}$ which results $R^{\alpha^{\prime}} = (h^{\alpha^{\prime}\alpha}) R^{\alpha} (h^{\alpha^{\prime}\alpha})^{-1}$. Therefore, to define global forms, it is sufficient to consider polynomials invariant under conjugation by general linear group. Hence, because of the multiplicative property of Berezinian(superdeterminant), one may consider $Ber(I+zR^{\alpha})$($z$ is a complex variable) which has a power expansion as follows(\cite{Voronov}):
\begin{equation}\label{power expansion}
Ber(I+zR^{\alpha}) = \sum_{k=0}^{\infty} c_{k}(R^{\alpha})z^{k}
\end{equation}
where $c_{k}(R^{\alpha}) = Tr \wedge^{k} R^{\alpha}$ and by $Tr$ we mean the supertrace.
\\
Now, it is necessary to show that $Tr \wedge^{k} R^{\alpha}$ is closed. To this end, note that we have
\begin{equation*}
Ber(I+zR^{\alpha}) = \exp\big(Tr \ln(I+zR^{\alpha})\big) = \exp \big(zTr(R^{\alpha}) - \dfrac{z^{2}}{2} Tr((R^{\alpha})^{2})+ \dfrac{z^{3}}{3} Tr((R^{\alpha})^{3}) -\cdots\big).
\end{equation*}
A comparison with \eqref{power expansion} shows that the coefficients $c_{k}(R^{\alpha})$ can be expressed as polynomials in $s_{k}(R^{\alpha}):= Tr (R^{\alpha})^{k}$. By \cite{Voronov}, $c_{k}(R^{\alpha})$ and $s_{k}(R^{\alpha})$ are connected by the following relation:
\begin{equation*}
c_{k+1} = \dfrac{1}{k+1} (s_{1}c_{k} - s_{2}c_{k-1} + \cdots + (-1)^{k}s_{k+1})
\end{equation*}
Thus, we should show that $Tr(R^{\alpha})^{k}$ is closed. To this end, we prove the following proposition:
\begin{proposition}
Let $\omega^{\alpha}$ and $R^{\alpha}$ be supermatrices defined by \eqref{omega def} and \eqref{Omega defi} respectively. We have
\begin{equation*}
dR^{\alpha} = \omega^{\alpha} \wedge R^{\alpha} - R^{\alpha} \wedge \omega^{\alpha}= [\omega^{\alpha}, R^{\alpha}]
\end{equation*}
\end{proposition}
\textit{Proof}:
Using the formula $R^{\alpha} = d\omega^{\alpha} - \omega^{\alpha} \wedge \omega^{\alpha}$, we get
\begin{align*}
dR^{\alpha} &= dd\omega^{\alpha} - d\omega^{\alpha} \wedge \omega^{\alpha} +\omega^{\alpha} \wedge d\omega^{\alpha} = -(R^{\alpha} + \omega^{\alpha} \wedge \omega^{\alpha})\wedge \omega^{\alpha} + \omega^{\alpha} \wedge (R^{\alpha} +\omega^{\alpha} \wedge \omega^{\alpha}) 
\\
                             &=-R^{\alpha} \wedge \omega^{\alpha} + \omega^{\alpha} \wedge R^{\alpha} = [\omega^{\alpha}, R^{\alpha}].
\end{align*}
\begin{flushright}
$\square$
\end{flushright}
Hence, we have
\begin{equation*}
dTr\big((R^{\alpha})^{k}\big) = \sum_{i+j=k-1} Tr\big((R^{\alpha})^{i} (dR^{\alpha}) (R^{\alpha})^{j}\big) = \sum_{i+j=k-1} Tr\big((R^{\alpha})^{i} [\omega^{\alpha}, R^{\alpha}] (R^{\alpha})^{j}\big) = Tr  [\omega^{\alpha}, (R^{\alpha})^{k}] = 0.
\end{equation*}

One may consider $c_k(R^\alpha)$ as a de Rham representative of the $k$-th $\nu$-class of the canonical super vector bundle $^{\nu}\gamma_{k|l}$. Since de Rham cohomology has finite dimension, there exists some index $i$ such that the classes of $c_k(R^\alpha)$, (for $k>i$), are equal to zero.
\section{Conclusion}
In this paper, via an analytic approach, the first $\nu$-class of a canonical $1|0$-super line bundle over $\nu$-projective space is represented by a 2-cocycle as an element of generalized de Rham cohomology. At the end, in a similar way of representation of Chern classes in common geometry, by considering canonical super vector bundles over $\nu$-grassmannian, the closed forms $c_k(R^\alpha)$ are introduced which may be considered as representatives of $\nu$-classes of rank higher than 1. In common geometry, for each Chern class, there is a representative in terms of Schubert cycles. It seems that a similar phenomenon may occur in supergeometry in the case of finding a proper generalization of Schubert cycles.
\addcontentsline{toc}{chapter}{References \hspace{13cm}}

{\small
\textit{E-mail addresses}: \href{mailto: m.roshandel@iasbs.ac.ir}%
{m.roshandel@iasbs.ac.ir} (M. Roshandelbana),
\\ 
     \hspace*{2.5cm}                                 \href{mailto: varsaie@iasbs.ac.ir}%
{varsaie@iasbs.ac.ir} (S. Varsaie).


\begin{thebibliography}{100}
\small
\bibitem{Afshari}
Afshari M. J. and Varsaie S., \textit{Gauss Supermaps and Their Homotopies}(Submitted).

\bibitem{Bahadory}
Bahadorykhalili F., Mohammadi M. and Varsaie S., \textit{A Class of Homogeneous Superspaces Associated
to Odd Involutions}, Periodica Mathematica Hungerica, (Accepted Feb 2020) .

\bibitem{Bartocci}
Bartocci C. and Bruzzo U., \textit{Cohomology of The Structure Sheaf of Real and Complex Supermanifolds}, Journal of Mathematical Physics \textbf{29}, 1789-1794 (1988).

\bibitem{Bartocci book}
Bartocci C., Bruzzo U. and Hernandez Ruiperez D., \textit{The Geometry of Supermanifolds}, Kluwer Academic Publisher, 1991.

\bibitem{Bruzzo}
Bruzzo U. and Hernandez Ruiperez D., \textit{Characteristic Classes of Super Vector Bundles}, Journal of Mathematical Physics \textbf{30}, 1233-1237 (1989).

\bibitem{DeWitt}
De Witt B., \textit{Supermanifolds}, Cambridge Univ. Press, London, 1984.

\bibitem{Griffith}
Griffith P. and Harris J., \textit{Principles of Algebraic Geometry}, Pure and Applied Mathematics, Wiley-Interscience, New York, 1978.

\bibitem{Hirzebruch}
Hirzebruch F., \textit{Topological Methods in Algebraic Geometry}, Springer Berlin Heidelberg, 1978.

\bibitem{Husemoller}
Husemoller D., \textit{Fibre bundles}, Springer-Verlag, 1994.

\bibitem{Voronov}
Khudaverdian H.M. and Voronov T. T., \textit{Berezinians, Exterior Powers and Recurrent Sequences}, Lett. Math. Phys., Vol. 74, 201-228 (2005).
 
\bibitem{Kostant}
Kostant B., \textit{Graded Manifolds, Graded Lie Theory, and Prequantization}, Lect. Notes Math., Vol. 570, 177-306 (1977). 

\bibitem{Landi}
Landi G., \textit{Projective Modules of Finite Type Over The Supersphere $S^{2,2}$}, Differential Geometry and its Application \textbf{14}, 95-111 (2001).

\bibitem{Leites}
Leites D. A., \textit{Introduction To The Theory of Supermanifolds}, Russian Mathematical Surveys, \textbf{35}, No. 1, 1-64 (1980).

\bibitem{Manin}
Manin Yu. I., \textit{Gauge Field Theory and Complex Geometry}, Springer-Verlag, 1984.
 
\bibitem{Manin and Penkov}
Manin Yu. I. and Penkov I. B., \textit{The Formalism of Left and Right Connections on Supermanifolds}, Lectures on Supermanifolds, Geometrical Methods of Conformal Groups, World Scientific (1989).

\bibitem{Quillen}
Quillen D., \textit{Superconnections and The Chern Character}, Topology, \textbf{24}, 89-95 (1985).

\bibitem{Rogers}
Rogers A., \textit{A Global Theory of Supermanifolds}, Journal of Mathematical Physics, \textbf{21}, 1352-1365 (1980).

\bibitem{Taniguchi}
Taniguchi T., \textit{ADHM Construction of Super Yang-Mills Instantons}, Journal of Geometry and Physics \textbf{59}, 1199-1209 (2009).

\bibitem{Varadarajan}
Varadarajan V. S., \textit{Supersymmetry For Mathematician: An Introduction}, American Mathematical Society, 2004.

\bibitem{Dr2}
Varsaie S., \textit{Homotopy Classification of Super Line Bundles}, XXXI Workshop on Geometric Methods in Physics, At Bialowieza, Poland, (2012).

\bibitem{Dr5}
Varsaie, S., \textit{On A Theory of Chern Classes For Super Vector Bundles}, Int. J. Geom. Methods Mod. Phys., \textbf{5}, No. 3, 287-295 (2008).

\bibitem{Dr1}
Varsaie S., \textit{$\nu$-Classes}, Int. J. Geom. Methods Mod. Phys., \textbf{9}, No. 4 (2012).

 \bibitem{Voronov-Manin-Penkov}
Voronov A. A., Manin Yu. I. and Penkov I. B., \textit{Elements of Supergeometry}, J. Math. Sci. \textbf{51}, 2069-2083 (1990).

\end{thebibliography}
\end{document}